\documentclass[preprint]{elsarticle}

\usepackage{hyperref}

\usepackage[utf8]{inputenc}  
\usepackage{amsmath,amsfonts,amssymb}
\usepackage{bm} 
\usepackage{tabularx} 
\usepackage{booktabs} 
\usepackage{mathtools} 
\usepackage{algorithm,algpseudocode}  
\usepackage[caption=false]{subfig} 

\usepackage{tikz}
\usetikzlibrary{shapes,arrows,calc,patterns}

\usepackage[a4paper,tmargin=2.5cm,bmargin=2.5cm,lmargin=2.5cm,rmargin=2.5cm]{geometry}

\newdefinition{remark}{Remark}

\newcommand{\domain}{D}
\newcommand{\transp}{\textsf{T}}

\newcommand{\relu}{\sigma}

\newcommand{\nodes}{N}

\newcommand{\wellnodes}{N^{\text{w}}}
\newcommand{\mannodes}{N^{\text{m}}}
\newcommand{\sepnodes}{N^{\text{s}}}

\newcommand{\edges}{E}
\newcommand{\dedges}{E^{\text{d}}}
\newcommand{\risedges}{E^{\text{r}}}

\newcommand{\oil}{\text{oil}}
\newcommand{\gas}{\text{gas}}
\newcommand{\gor}{\text{gor}}
\newcommand{\wor}{\text{wor}}
\newcommand{\water}{\text{wat}}

\newcommand{\phases}{C}

\begin{document}

\begin{frontmatter}

\title{ReLU networks as surrogate models in mixed-integer linear programs}

\author[address1]{Bjarne Grimstad\corref{correspondingauthor}}
\ead{bjarne.grimstad@gmail.com}
\author[address2]{Henrik Andersson}

\address[address1]{Solution Seeker AS, Gaustadalléen 21, 0349 Oslo, Norway}
\address[address2]{Department of Industrial Economics and Technology Management, Norwegian University of Science and Technology, NO-7491 Trondheim, Norway}

\cortext[correspondingauthor]{Corresponding author}


%

\begin{abstract}
We consider the embedding of piecewise-linear deep neural networks (ReLU networks) as surrogate models in mixed-integer linear programming (MILP) problems. A MILP formulation of ReLU networks has recently been applied by many authors to probe for various model properties subject to input bounds. The formulation is obtained by programming each ReLU operator with a binary variable and applying the big-M method. The efficiency of the formulation hinges on the tightness of the bounds defined by the big-M values. When ReLU networks are embedded in a larger optimization problem, the presence of output bounds can be exploited in bound tightening. To this end, we devise and study several bound tightening procedures that consider both input and output bounds. Our numerical results show that bound tightening may reduce solution times considerably, and that small-sized ReLU networks are suitable as surrogate models in mixed-integer linear programs. 
\end{abstract}


\end{frontmatter}


\section{Introduction}
Access to larger datasets and more computational power have created vast possibilities for the application of machine learning algorithms to data-driven decision support systems studied in operations research and mathematical programming. We are now able to exploit the data and build highly complex models to describe nonlinear phenomena. A challenge in mathematical optimization is how to use these models in an optimization framework and how to efficiently solve the problems that arise. 

Models that act as substitutes for nonlinear relationships in mathematical programming problems are often called \emph{surrogate models}. The study of surrogate models is important for several reasons: 1) it allows optimization of functions with an unknown/hidden mathematical structure, for which derivative information is unavailable; 2) it may reduce solution times when complex functions can be substituted by simpler surrogate models with better properties for optimization (e.g. smoothness and linearity); and 3) for problems with unknown functions that are in some sense expensive to sample, surrogate models can be used to handle the trade-off between exploration and exploitation \cite{Jones1998,Ghavamzadeh2015}. The view of surrogate modeling is thus useful when considering the challenge of embedding machine learning models in mathematical programs.

Functions with the above properties, which are accessible only via sampling, are called ``black-box'' functions. Optimization problems containing a mix of black-box functions and explicitly known functions (for example algebraic expressions) are often referred to as ``grey-box'' problems \cite{Beykal2018}. Problems with this setup are commonly encountered in process optimization. A promising approach to the optimization of grey-box problems is to approximate and replace the black-box functions with surrogate models. The approximate problem can then be efficiently solved by derivative-based optimization methods. This is the approach taken by model-based derivative-free optimization (DFO) methods \cite{Conn2009}.

Model-based DFO methods can be categorized as \emph{local} and \emph{global} DFO methods \cite{Bhosekar2018}. Local DFO methods are specialized to search efficiently for a local optimum and often rely on building rather simple surrogate models. In contrast, global DFO methods build surrogate models over the entire feasible region, and attempt to locate a global optimum. As global surrogate models are necessarily required to be non-convex, global optimization algorithms are usually employed to solve the problem \cite{Boukouvala2017}. Most DFO methods incorporate a sampling algorithm that decides at which points to sample the black-box functions. An initial sample is used to fit the surrogate models, and subsequent samples can be taken to locally refine the surrogate models before reoptimizing the problem. A review of surrogate modeling and optimization can be found in \cite{Rios2013}, and more recently in \cite{Bhosekar2018}. \cite{Bhosekar2018} also provide an overview of software for building various surrogate models, including some of the surrogate models we mention below.

Many of the classic function approximation methods have been applied to build surrogate models for process analysis and optimization \cite{Bhosekar2018}. Examples include radial basis functions \cite{Regis2007}, support vector machines \cite{Ciccazzo2014}, tensor-product splines \cite{Grimstad2016b}, and artificial neural networks \cite{Fahmi2012, Ye2019}. Generally, the resulting surrogate models are constructed to be continuous and smooth to facilitate optimization by nonlinear (global) solvers. An important class of surrogate models is the piecewise-linear (PWL) models, which naturally lend themselves to mixed-integer linear programming (MILP) formulations. Various MILP formulations exist for PWL models of one or multiple variables, including the ones listed in \cite{Vielma2010}. The starting point of these formulations is a set of sample points to be interpolated. As will be discussed subsequently, some of these formulations put restrictions on the sample structure in addition to the interpolation constraints. This limits their applicability to low-dimensional functions and relatively small datasets with a few thousand data points.


When modeling on large and noisy datasets, for example to fit global surrogate models, regression methods are more suitable than interpolation methods. The modeler can then leverage tools for statistical analysis and machine learning to build accurate PWL approximations by controlling model complexity to avoid overfitting the data. With appropriate formulations, the resulting regression models can be brought into a MILP framework for optimization. In this spirit, many works have been published recently on how to embed trained machine learning models in optimization \cite{Martinez2017, Elmachtoub2017, Misic2017, Biggs2018}. For example, \cite{Martinez2017} recently employed a MILP formulation for multivariate adaptive regression splines (MARS), while \cite{Elmachtoub2017} developed a framework for parameter prediction that incorporates the objective and constraints of the optimization problem explicitly. The desire to bring surrogate models into a MILP framework stems from the demonstrated effectiveness of MILP solution methods in solving large-scale instances of complex problems; see for example \cite{Grossmann2012,Andersson2015b,Bach2016,Veenstra2017}.


We follow a similar path in this work by employing a MILP formulation for the general class of \emph{ReLU networks}. These are networks composed of max-affine spline operators (piecewise-linear and convex operators), which include: affine operators, ReLU-type activations, absolute value activations, convolution operators, and max/mean pooling. By construction, ReLU networks are affine spline operators, which are piecewise linear, but not necessarily convex \citep{Balestriero2018}. In this class of deep neural networks (DNNs) we find widely applied architectures such as convolutional neural networks \cite{Goodfellow2016}, residual neural networks \cite{He2016}, inception networks \cite{Szegedy2016}, maxout networks \cite{Goodfellow2013}, network-in-networks \cite{Lin2013}, and their variants using max-affine spline operators. These networks are used to model complex relationships and often achieve state-of-the-art performance in terms of modeling accuracy \cite{Glorot2011, He2015, Schmidhuber2015}. In this work we consider ReLU networks composed solely of affine operators and ReLU activations, since these operators are sufficient for many function approximation tasks; a shallow ReLU network with one hidden layer is a universal approximator \cite{Sonoda2017}.

An exact MILP formulation of a ReLU network can be obtained by programming each ReLU operator with a binary variable and applying the big-M method. This formulation has recently been applied to formal verification \cite{Bunel2017, Cheng2017, Dutta2017, Tjeng2017, Fischetti2018}, to count linear regions \cite{Serra2018}, and to compress DNNs \cite{Kumar2019}. Common for these works is consideration of a single ReLU network subject to input bounds. These studies show that the application of bound tightening techniques to compute big-M values can reduce solution times considerably. 

We use the same formulation to solve optimization problems with multiple ReLU networks embedded. The problems we consider can be viewed as grey-box problems, where ReLU networks are used as global surrogate models. In this setting, multiple ReLU networks are interlinked by algebraic constraints, eliciting bounds on the network outputs. The studies listed above are limited to bound tightening techniques that propagate input bounds forward through a single network. To take advantage of the output constraints present in grey-box problems, we study stronger bound tightening techniques that are capable of propagating output bounds backwards. The consideration of bound tightening in the presence of output bounds, and the optimization of problems with multiple ReLU networks, distinguishes our study from previous works.

The focus of this paper is to show how ReLU networks can be used to model complex phenomena in a mixed-integer linear framework and how bound tightening affect the performance of the big-M formulation. Our main contributions are summarized as follows:
\begin{itemize}
	\item[--] A discussion about the advantages of using ReLU networks to model complex phenomena in a MILP framework.
	\item[--] A framework for tightening bounds that unify different bound tightening procedures for output-bounded ReLU networks under one general procedure.
	\item[--] A computational study showing the strength of the bound tightening procedures, but also the important trade-off between time spent on bound tightening and the solution time of the optimization problem.
	\item[--] A real application where bound tightening makes the difference between not finding a feasible solution and solving the problem to optimality within the practical limit on solution time.
\end{itemize}

The literature contains numerous applications of neural networks to model and locally optimize complex processes; a list all previous works would be quite comprehensive so we only point to a few examples here \cite{Fernandes2006,Santanna2017,AL-Qutami2018}. We firmly believe that ReLU networks, which are the de-facto standard in deep learning \cite{Goodfellow2016}, is an important class of surrogate models that deserves being studied. Our work is also motivated by the desire to pair state-of-the-art software for machine learning with mixed-integer linear optimization. In particular, we wish to leverage powerful modeling frameworks like TensorFlow \cite{tensorflow2015-whitepaper} to build PWL surrogate models, and then incorporate these in a MILP model that we can solve with a state-of-the-art MILP solver. 

Our assessment is on the applicability of ReLU networks as global surrogate models in a grey-box approach to process optimization. We do not consider sampling to refine and reoptimize the problem, but focus on building highly accurate approximations to be used in a single global optimization. However, our findings can be used to develop more sophisticated model-based global DFO methods based on ReLU networks.

\subsection{Structure of the paper}
We begin in Section \ref{sec:pwl-approximation} by discussing the various aspects of building piecewise-linear approximations of nonlinear functions, and embedding these as surrogate models in optimization. We highlight the advantages of modeling multi-variable functions via ReLU networks, as opposed to using interpolation on simplices, which has been the traditional approach in MILP literature. After presenting the MILP formulation in Section \ref{sec:dnn-milp-model}, we devise several optimization-based bound tightening procedures for output-bounded ReLU networks in Section \ref{sec:bound-tightening}. To study the procedures' efficiency, and the feasibility of using ReLU networks as surrogate models in process optimization, we present a computational study in Section \ref{sec:numerical-study}. Our test suite includes a challenging oil production optimization problem, for which the complicated physics of multiphase flow in pipes is modeled by ten ReLU networks. Concluding remarks and promising research directions are given in Section \ref{sec:concluding-remarks}.

\section{Piecewise-linear approximations of nonlinear functions}
\label{sec:pwl-approximation}
We consider the modeling of a nonlinear function $f : \domain \subset \mathbb{R}^n \to \mathbb{R}$, on a compact domain $\domain \subset \mathbb{R}^n$ (we assume that $\domain$ is a polytope). To incorporate $f$ in a MILP model, it must be approximated by a piecewise-linear function. For cases where $f$ is unknown, nonseparable (for $n \geq 2)$, or highly complex, a general approach to the approximation of $f$ is to sample it on the domain $\domain$, and then build a piecewise-linear approximation from the sample. In cases where $f$ is not a mathematical construct, for example if $f$ is a real process, the sample may represent a set of experiments. In general, a sample consists of non-structured sample points scattered in $\domain$.

A common restriction of many approximation methods is that $f$ must be sampled on a rectilinear grid $G$, covering a hyper-rectangular domain $\domain$. A rectilinear grid on which variable $x_i$ is partitioned into $m_i$ intervals, results in $\prod_{i=1}^{n} m_i$ boxes (hyperrectangles). Clearly, the number of boxes grows exponentially with $n$. To obtain a piecewise linear model that interpolates the sample points, each of these boxes must be divided into a set of simplices. The minimum number of simplices needed to triangulate an $n$-dimension hypercube is 1, 2, 5, 16, 67, 308, and 1493 for $1\leq n \leq 7$, respectively, see \cite{Hughes1996}. In general, an upper bound is given by $n!$. Thus, a partitioning of $\domain$ may consist of up to $n! \prod_{i=1}^{n} m_i$ simplices.

Table \ref{tab:number-of-simplices} lists the lower bound on the number of simplices resulting from a rectilinear grid partitioning. The bounds are given for $1 \leq n \leq 7$, with each variable partitioned into ten intervals. Reading the table, we see that the exponential increase in the number of boxes (or simplices) prohibits any practical use of such a partitioning for $n > 3$. For example, \cite{Misener2010} and \cite{Vielma2010} consider grids for $n \leq 3$.

\begin{table}[ht]
\footnotesize
\caption{Number of simplices resulting from a rectilinear grid partitioning of the domain, where each variable is partitioned into ten intervals. We have used the lower bound on the number of simplices per box.}
\begin{tabularx}{\textwidth}{l*3{>{\raggedleft\arraybackslash}X}}
\toprule
Dimension $n$ & \# boxes & \# simplices/box & \# simplices \\
\midrule 
1 & 10 & 1 & 10 \\
2 & 100 & 2 & 200 \\
3 & 1 000 & 5 & 5 000 \\
4 & 10 000 & 16 & 160 000 \\
5 & 100 000 & 67 & 6 700 000 \\
6 & 1 000 000 & 308 & 308 000 000 \\
7 & 10 000 000 & 1 493 & 14 930 000 000 \\
\bottomrule \noalign{\smallskip}
\end{tabularx}
\label{tab:number-of-simplices}
\end{table}

Another disadvantage with rectilinear sampling, stemming from its structure, is its inability to \emph{locally} control the sampling resolution. That is, increasing the sampling resolution in a region of $\domain$ where $f$ has important features (e.g. large gradients), will also increase it in other regions where it may not be necessary to sample densely. Rectilinear sampling is also subject to high sparsity of sample points in higher dimensions. Although this aspect of the ``curse of dimensionality'' applies to all sampling methods, it is becomes especially prominent for rectilinear sampling due to its inability to locally control the sampling resolution. 


For high-dimensional problems, there exist sampling algorithms with better space-filling properties than rectilinear sampling, such as Latin hypercube sampling. The sampling strategy can either be static, as in rectilinear sampling, or sequential. With sequential sampling, sample points are successively selected to optimize a space-filling criterion or to minimize the (expected) approximation error of a surrogate model \cite{Crombecq2011, Gorissen2010, Herten2015}. In the latter case, the properties of the surrogate model can be exploited to reduce the sample size \cite{Kieslich2018}. These sampling stategies cannot escape the ``curse of dimensionality'', but may scale more graciously to higher dimensions by being economic in the placement of the sample points. This ability is especially important when $f$ is expensive to evaluate and there are practical limits to the sample size. We refer the reader to \cite{Garud2017} for a thorough review of modern sampling algorithms.

Sampling algorithms like those described above typically result in sample points that are scattered in $\domain$ without any particular structure. The same holds true for most samples obtained from running real experiments. A piecewise-linear model can be obtained from such samples via triangulation of the sample points, resulting in a simplex partitioning of $\domain$. There exist several MILP models for PWL functions that have no requirements on the family of polytopes/simplices produced by the triangulation, and which can be used with a triangulation scheme for scattered sample points. One example is the DLog model in \cite{Vielma2010}, which uses a logarithmic number of binary variables to model a PWL function. This model has shown good computational performance for higher-dimensional functions ($n \geq 2$), compared to MILP models with a linear number of binary variables. However, as is shown subsequently, this modeling approach is not feasible for the oil production optimization case considered in this paper.

There exist MILP formulations for PWLs that achieve higher efficiency by putting requirements on the partitioning scheme. E.g., the Log model of \cite{Vielma2010} requires a triangulation scheme compatible with the $J_1$ (``Union Jack'') triangulation, while SOS2 models require a rectilinear partitioning \cite{Beale1970}. In general, these highly structured partitioning schemes do not conform with a sample of scattered points and are also subject to the dimensionality issues discussed above.

So far we have only discussed approximation methods that interpolate the sample points.\footnote{We note here that for $n \geq 2$, SOS2 models do not interpolate all sample points.} If this requirement is relaxed, we may consider regression methods using piecewise-linear splines. These methods solve a regression problem by minimizing the approximation error on the sample. We may divide these methods into two classes: non-adaptive methods for which the domain partitioning is fixed (given by break-points or knots); and adaptive methods where the model parameters define the polytopic partition of the domain and the linear relationship in each polytope. MARS and ReLU networks are examples of adaptive methods. 

Regression models offer a nice alternative to the interpolation methods in that their capacity (number of linear pieces) is easily controlled. This is a key property when building predictive models, as it allows the modeler to strike a balance between underfitting and overfitting the data. Additionally, it can be advantageously used to reduce \emph{unnecessary} complexity of surrogate models in a mathematical optimization problem. As the number of binary variables in a MILP formulation scales directly with the number of partition regions (number of linear pieces), it is reasonable to expect that surrogate models with fewer regions can be optimized faster. In accordance with the principle of \emph{Occam's razor}, unnecessary complexity, which here translates to regions that do not significantly lower the model accuracy, is favorably avoided in both the modeling and optimization setting.

\subsection{ReLU networks as piecewise-linear surrogate models}
Below, we highlight some important properties of ReLU networks that make them attractive as surrogate models in optimization.

\paragraph{Piecewise-linear and continuous} 
ReLU networks are composed of max-affine spline operators and are piecewise-linear and continuous functions \cite{Balestriero2018}. This means that a ReLU network can be \emph{exactly} formulated and optimized in a MILP framework. Networks composed of other nonlinear operators, such as sigmoid or hyperbolic functions, can only be formulated approximately in a MILP framework via piecewise-linear approximation of the operators.

\paragraph{Suited for modeling of high-dimensional and complex functions} 
It is easy to control the capacity of a ReLU network by adjusting its depth and width. As depth and width increases, so does the number of nodes (also called neurons) and effectively the capacity of the model. An upper bound on the number of regions (linear pieces) generated by a neural network with $L$ hidden layers of width $W$ and inputs $x \in \mathbb{R}^{n_0}$, is $O(W^{L n_0})$ \cite{Raghu2016}. For fixed $W$ and $n_0$, this upper bound indicates that the number of regions grows exponentially in the network depth $L$. The bound shows that deep networks can be more expressive than shallow networks (can be seen by keeping $WL$ constant and varying $L$). This observation has been used to partly explain the success of DNNs \cite{Serra2018}.

\begin{figure}[ht]
    \centering
    \includegraphics[width=1.0\textwidth]{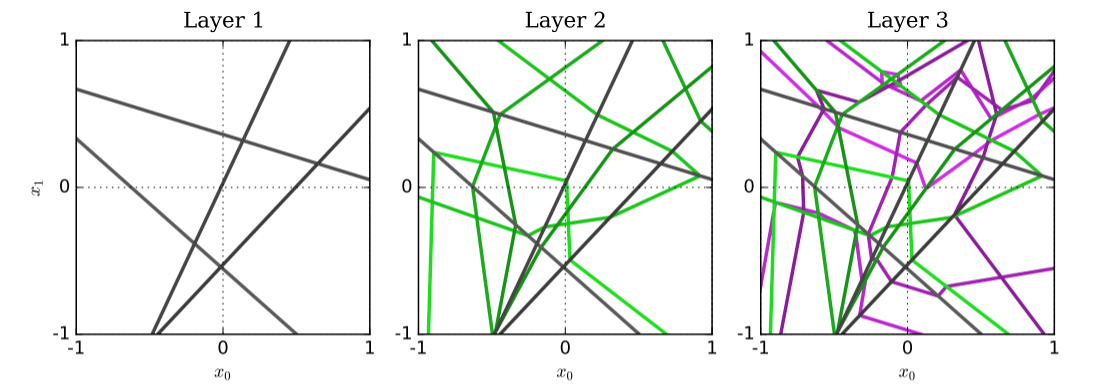}
	\caption{Example of domain partition generated by a multilayered neural network; adopted from \cite{Raghu2016}.}
	\label{fig:nn-regions}
\end{figure}

\paragraph{Adaptive partitioning of the domain} 
Changes in the parameters of a ReLU network automatically induce changes in the partition of the domain \citep{Balestriero2018}. When the network parameters are updated during training, regions of the domain where the data has important features are automatically described with smaller polytopes, while other parts of the domain are described by larger polytopes. Figure \ref{fig:nn-regions} shows an example of a domain partition, which is refined from one layer to the next.

\paragraph{Can be trained on scattered and noisy data} 
The process of training a ReLU network does not put requirements on the distribution of the data points and is thus compatible with any sampling algorithm. Furthermore, an arsenal of regularization techniques is available for neural networks in general. These techniques allow neural networks to generalize from noisy data \cite{Goodfellow2016}. Scattered and noisy data are issues that occur in many practical applications. 

\paragraph{Scale to large datasets} Neural networks are usually trained with (mini-batch) stochastic gradient descent algorithms \cite{Bottou2018}. These algorithms optimize over the parameters $\Theta$ of a neural network $f_{\Theta}$ in search of an approximation of some true function $f$. For regression tasks, it is common to minimize the empirical loss
\begin{equation*}
\mathcal{L} = \frac{1}{N} \sum\limits_{i = 1}^{N} \Vert y^{(i)} - f_{\Theta}(x^{(i)}) \Vert_{2}^{2} + \lambda \Vert \Theta \Vert_{2}^{2},
\end{equation*}
where $\{x^{(i)}, y^{(i)}\}_{i=1}^{N}$ are $N$ samples of $f$. The first term specifies the mean squared error, while the second is an $L_2$ regularization term that penalizes large values of $\Theta$. In addition to the network width and depth, the hyperparameter $\lambda \in \mathbb{R}_{\geq 0}$ is used to control model complexity and encourage a low generalization error.

During optimization, the parameters are updated using gradients $\partial \mathcal{L}/\partial \Theta$ computed using backward-propagation. The batch size $N$ can be selected to control the computation and memory resources required to evaluate $\mathcal{L}$ and its gradients. This enables training of deep neural networks on large datasets with millions of data points. 

\paragraph{Many software tools available for modeling} 
A plethora of software tools exists for designing, training and performing inference tasks with ReLU networks. In this work we use the machine learning framework TensorFlow \cite{tensorflow2015-whitepaper}.

\section{A 0--1 MILP model for ReLU networks}
\label{sec:dnn-milp-model}
We consider a ReLU network with $K+1$ layers, numbered from 0 to $K$. Layer 0 is the input layer (usually not counted as a layer), while the last layer $K$ is the output layer. Each layer $k \in \{0, \ldots, K\}$ has $n_k$ nodes, numbered from $1$ to $n_k$. For simplicity we consider networks of fully connected layers with parameters $\theta^{k} = (W^k, b^k)$ for $k \in \{1, \ldots, K\}$, where $W^{k} \in \mathbb{R}^{n_{k} \times n_{k-1}}$ and $b^{k} \in \mathbb{R}^{n_k}$. We state the network as $f_{\Theta} : \mathbb{R}^{n_0} \to \mathbb{R}^{n_K}$, where $K$, $n_0$, and $n_K$ are implied by the set of parameters $\Theta := \left\{ \theta^1, \ldots, \theta^K \right\}$.

Let $x^{k} \in \mathbb{R}^{n_k}$ be the output of layer $k$, and $x_{j}^{k}$ the output of the $j$-th node for $j=1,\ldots,n_k$. Following this notation, $x_{j}^{0}$ is the $j$-th input value, and $x_{j}^{K}$ is the $j$-th output value (out of $n_K$) of the network. For each (hidden) layer $k \in \{1,\ldots,K-1\}$, the output vector $x^k$ is computed as
\begin{equation}
x^k = \sigma(W^k x^{k-1} + b^k),
\label{eq:ReLU-layer}
\end{equation}
where we denote by $\relu(y) := \max \{0, y\}$ (componentwise) the ReLU activation function for a real vector $y$. The output of the network, $x^{K} \in \mathbb{R}^{n_K}$, is given by the affine equation $x^K = W^K x^{K-1} + b^K$.

\begin{figure}[ht]
	\begin{centering}
		\begin{tikzpicture}[
		cir/.style={circle,inner sep=0pt,minimum size=10mm,draw=black},
		dot/.style={circle,inner sep=0pt,minimum size=1mm,fill=black},
		rec/.style={draw,minimum width=\wid cm,minimum height=0.55cm}
		]
		\node at (0,2.8) {Input layer};
		\node at (0,2.3) {$0$};
		\node at (2,2.3) {$1$};
		\node at (4.5,2.3) {$k-1$};
		\node at (6.5,2.3) {$k$};
		\node at (9,2.3) {$K-1$};
		\node at (11,2.3) {$K$};
		\node at (11,2.8) {Output layer};
		\draw [|-, dashed] (2,2.8) to (4.0,2.8);
		\node at (5.5,2.8) {Hidden layers};
		\draw [-|, dashed] (7,2.8) to (9,2.8);	

		\node[cir] (n0) at (0,1.5) {$x^0$};
		\node[cir] (n1) at (2,1.5) {$x^1$};
		\node (na) at (2.8,1.5) {};
		\node (nb) at (2.9,1.5) {};
		\node[cir] (ni) at (4.5,1.5) {$x^{k-1}$};
		\node[cir] (nj) at (6.5,1.5) {$x^k$};
		\node (ne) at (7.8,1.5) {};
		\node (nf) at (7.9,1.5) {};
		\node[cir] (nk) at (9,1.5) {$x^{K-1}$};
		\node[cir] (nK) at (11,1.5) {$x^K$};
		
		\draw [->] (n0) to (n1);
		\draw [-,dashed] (n1) to (na);
		\draw [->,dashed] (nb) to (ni);
		\draw [->] (ni) to (nj);
		\draw [-,dashed] (nj) to (ne);
		\draw [->,dashed] (ne) to (nk);
		\draw [->] (nk) to (nK);

		\draw[dotted] (3.2,0.9) rectangle (7.4,2.1);
		\draw[dotted] (2.2,-6.2) rectangle (8.9,0.0);
		\draw[dotted] (3.2,0.9) -- (2.2,-0.0);
		\draw[dotted] (7.4,0.9) -- (8.9,-0.0);

		\node at (1.6,1.5) [anchor=north east] {$\theta^1$};
		\node at (4.1,1.5) [anchor=north east] {$\theta^{k-1}$};
		\node at (6.1,1.5) [anchor=north east] {$\theta^k$};
		\node at (8.6,1.5) [anchor=north east] {$\theta^{K-1}$};
		\node at (10.6,1.5) [anchor=north east] {$\theta^K$};

		\node[cir] (ni1) at (3,-1) {$x^{k-1}_1$};
		\node[cir] (nj1) at (8,-1) {$x^{k}_1$};
		\node[cir] (ni2) at (3,-3) {$x^{k-1}_2$};
		\node[cir] (nj2) at (8,-3) {$x^{k}_2$};
		\node[cir] (nin) at (3,-5.5) {$x^{k-1}_{n_{k-1}}$};
		\node[cir] (njn) at (8,-5.5) {$x^{k}_{n_k}$};

		\node[dot] at (3,-3.7) {};
		\node[dot] at (3,-4.0) {};
		\node[dot] at (3,-4.3) {};
		\node[dot] at (8,-3.7) {};
		\node[dot] at (8,-4.0) {};
		\node[dot] at (8,-4.3) {};

		\draw[->] (ni1) to node[very near end,above,sloped] {\small $W^k_{11}$} (nj1);
		\draw[->] (ni1) to node[very near end,above,sloped] {\small $W^k_{21}$} (nj2);
		\draw[->] (ni1) to node[very near end,above,sloped] {\small $W^k_{n_k1}$} (njn);
		\draw[->] (ni2) to (nj1);
		\draw[->] (ni2) to (nj2);
		\draw[->] (ni2) to (njn);
		\draw[->] (nin) to (nj1);
		\draw[->] (nin) to (nj2);
		\draw[->] (nin) to node[very near end,below,sloped] {\small $W^k_{n_k n_{k-1}}$} (njn);

		\node at (8,-0.6) [anchor=south] {\small $b^k_1$};
		\node at (8,-2.6) [anchor=south] {\small $b^k_2$};
		\node at (8,-5.1) [anchor=south] {\small $b^k_{n_k}$};

		\draw [-,dashed] (8.9,-1) to (nj1);
		\draw [-,dashed] (8.9,-1.2) to (nj1);
		\draw [-,dashed] (8.9,-1.4) to (nj1);
		\draw [-,dashed] (8.9,-2.8) to (nj2);
		\draw [-,dashed] (8.9,-3.0) to (nj2);
		\draw [-,dashed] (8.9,-3.2) to (nj2);
		\draw [-,dashed] (8.9,-5.1) to (njn);
		\draw [-,dashed] (8.9,-5.3) to (njn);
		\draw [-,dashed] (8.9,-5.5) to (njn);

		\draw [->,dashed] (2.2,-1) to (ni1);
		\draw [->,dashed] (2.2,-1.2) to (ni1);
		\draw [->,dashed] (2.2,-1.4) to (ni1);
		\draw [->,dashed] (2.2,-2.8) to (ni2);
		\draw [->,dashed] (2.2,-3.0) to (ni2);
		\draw [->,dashed] (2.2,-3.2) to (ni2);
		\draw [->,dashed] (2.2,-5.1) to (nin);
		\draw [->,dashed] (2.2,-5.3) to (nin);
		\draw [->,dashed] (2.2,-5.5) to (nin);

		\end{tikzpicture}
		\caption{The structure of a ReLU network. The top part shows the aggregated structure from the input vector $x^0$ to the output vector $x^K$. The lower part shows a disaggregated view of two layers where the variables and parameters of the network are shown.}
		\label{fig:ReLU-network}
	\end{centering}
\end{figure}
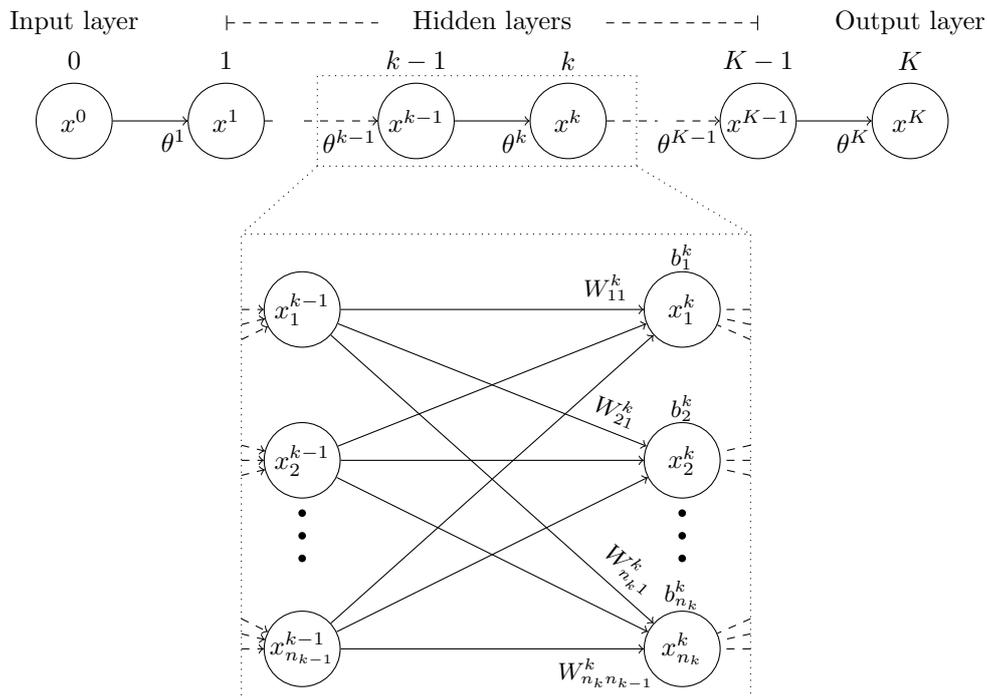

Figure \ref{fig:ReLU-network} shows the structure of a ReLU network. The output of layer $k$ is calculated using the input $x^{k-1}$ from the previous layer and the parameters $\theta^k = (W^k,b^k)$ according to \eqref{eq:ReLU-layer}.

The ReLU operator can be programmed in multiple ways. Following \citep{Fischetti2018}, we consider the linear equation
\begin{equation}
w^\transp y + b = x - s,
\end{equation}
where the output of the ReLU is decoupled into a positive part $x \geq 0$ and negative part $s \geq 0$. The output of the ReLU can then be obtained by imposing that at least one of the two terms $x$ and $s$ must be zero. Assuming that we may compute finite values $L$ and $U$ so that $L \leq w^\transp y + b \leq U$, we may explicitly program the ReLU logic via big-M constraints:
\begin{equation}
\begin{aligned} 
& x \leq Uz \\
& s \leq -L(1-z) \\
& z \in \{0, 1\}
\end{aligned}
\label{eq:bigm-constraints}
\end{equation}
where we have introduced a binary \emph{activation variable} $z$.


Using a binary activation variable for each node $(j,k)$, the network can be expressed with the following 0--1 MILP formulation:

\begin{subequations}
\begin{align}
& \text{Input layer} \notag \\
& L^{0} \leq x^{0} \leq U^{0}, \label{eq:nn-milp-input-bounds} \\
& \text{Hidden (ReLU) layers} \notag \\
& \begin{rcases} 
W^{k} x^{k-1} + b^{k} = x^{k} - s^{k} \\
x^{k}, s^{k} \geq 0 
\end{rcases} & k=1,\ldots,K-1, \label{eq:nn-milp-linear} \\
& z_{j}^{k} \in \{0, 1\} & k=1,\ldots,K-1, j=1,\ldots,n_k \label{eq:nn-milp-binary} \\
& \begin{rcases} 
x_{j}^{k} \leq U_{j}^{k} z_{j}^{k} \\
s_{j}^{k} \leq -L_{j}^{k} (1 - z_{j}^{k})
\end{rcases} & k=1,\ldots,K-1, j=1,\ldots,n_k \label{eq:nn-milp-bigm} \\
& \text{Output layer} \notag \\
& W^{K} x^{K-1} + b^{K} = x^{K}, \label{eq:nn-milp-output-layer} \\
& L^{K} \leq x^{K} \leq U^{K}. \label{eq:nn-milp-output-bounds}
\end{align}
\label{eq:nn-milp}
\end{subequations}

This is an exact formulation of the ReLU network, meaning that the output of this formulation always is the same as the output from the ReLU network for the same input. That is, for a fixed input vector $x^0$, all the other variables in \eqref{eq:nn-milp} are fixed, including the output variables $x^K$. The only ``degenerate'' solutions occur when the input to a node $(j,k)$ is zero, in which case the value of $z_{j}^{k}$ can arbitrarily be set to 0 or 1 without affecting the output.

There is one binary variable for each hidden node in the network, this means that the formulation scales linear with the number of hidden nodes. The constraint sets \eqref{eq:nn-milp-linear} and \eqref{eq:nn-milp-output-layer} can be written as $A_x x + A_s s = -b$, where $A_x$ and $A_s$ have a special block angular structure (depends on how we stack the variables $x$ and $s$), see Figure \ref{fig:ReLU-block-angular}. Unless sparsity is encouraged, for example by utilizing $L_1$ regularization, the parameters $W^k$ and $b^k$ are typically dense for a trained network.

\def\wid{1.25}

\begin{figure}[ht]
	\begin{centering}
	\begin{tikzpicture}[
		cir/.style={circle,inner sep=0pt,minimum size=1mm,fill=black},
		rec/.style={draw,minimum width=\wid cm,minimum height=0.55cm}
	]
	\node at (0*\wid,2.5) {$x^0$};
	\node at (1*\wid,2.5) {$s^1$};
	\node at (2*\wid,2.5) {$x^1$};
	\node at (3*\wid,2.5) {$s^2$};
	\node at (4*\wid,2.5) {$x^2$};
	\node at (6*\wid,2.5) {$x^{K-2}$};
	\node at (7*\wid,2.5) {$s^{K-1}$};
	\node at (8*\wid,2.5) {$x^{K-1}$};
	\node at (9*\wid,2.5) {$x^{K}$};
	\node[anchor=west] at (9.6*\wid,2.0) {$=-b^1$};
	\node[anchor=west] at (9.6*\wid,1.4) {$=-b^2$};
	\node[anchor=west] at (9.6*\wid,0.4) {$=-b^{K-1}$};
	\node[anchor=west] at (9.6*\wid,-0.2) {$=-b^{K}$};
	\node[rec] (rect) at (0,2) {$W^1$};
	\node[rec] (rect) at (\wid,2) {$I$};
	\node[rec] (rect) at (2*\wid,2) {$-I$};
	\node[rec] (rect) at (2*\wid,1.4) {$W^2$};
	\node[rec] (rect) at (3*\wid,1.4) {$I$};
	\node[rec] (rect) at (4*\wid,1.4) {$-I$};
	\node[cir] at (4.7*\wid,1.0) {};
	\node[cir] at (5.0*\wid,0.9) {};
	\node[cir] at (5.3*\wid,0.8) {};
	\node[rec] (rect) at (6*\wid,0.4) {$W^{K-1}$};
	\node[rec] (rect) at (7*\wid,0.4) {$I$};
	\node[rec] (rect) at (8*\wid,0.4) {$-I$};
	\node[rec] (rect) at (8*\wid,-0.2) {$W^{K}$};
	\node[rec] (rect) at (9*\wid,-0.2) {$-I$};
	\end{tikzpicture}
	\caption{The block-angular structure of the constraint set of a ReLU network.}
	\label{fig:ReLU-block-angular}
	\end{centering}
\end{figure}
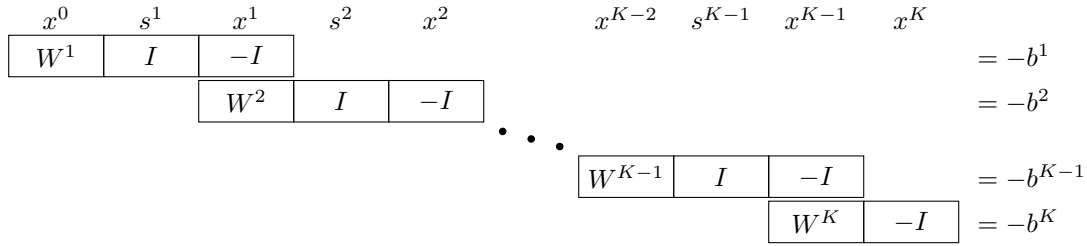

The tightness of an LP relaxation of \eqref{eq:nn-milp} depends heavily on the size of the big-M values $L^{k}$ and $U^{k}$, for layers $k=1,\ldots,K-1$. Weak LP relaxations result from large big-M values, and may severely hamper the efficiency of the solver. We explore several bound tightening procedures in the upcoming section. The computed bounds, $L_{j}^{k}$ and $U_{j}^{k}$, constrain the variables as follows:
\begin{equation}
\begin{aligned}
\max\{0, L_{j}^{k}\} &\leq x_{j}^{k} \leq \max\{0, U_{j}^{k}\} \\
\max\{0, -U_{j}^{k}\} &\leq s_{j}^{k} \leq \max\{0, -L_{j}^{k}\}
\end{aligned}
\label{eq:bound-update-rule}
\end{equation}

Before continuing, we remark that alternative formulations exist for the ReLU logic. The logic may for example be programmed via complementary constraints ($xs = 0$) or indicator constraints ($z = 0 \implies x \leq 0$ and $z = 1 \implies s \leq 0$). An advantage with these formulations is that they do not require explicit big-M values. However, these formulations tend to produce very hard mixed-integer instances that challenge state-of-the-art solvers \cite{Fischetti2018}.

\section{Bound tightening procedures}
\label{sec:bound-tightening}
We study procedures for computing valid bounds $B := \left\{ B^0, \ldots, B^K \right\}$, where $B^i := [L^i, U^i]$ for $i=0,\ldots,K$. We begin by defining the feasible set of a ReLU network $f_\Theta$ bounded by $B$ as
\begin{equation}
P(\Theta, B) := \left\{ (x, s, z) : \eqref{eq:nn-milp} \text{ is satisfied} \right\}.
\label{eq:feasible-set-milp}
\end{equation}
The projection of this set onto $(x^0, x^K)$, given by $F(\Theta, B) := \{ (x^0, x^K) : (x, y, z) \in P(\Theta, B) \}$, yields all possible input-output pairs of $f_\Theta$. The inputs $x^0$ and outputs $x^K$ are directly constrained by the bounding boxes $B^0$ and $B^K$, respectively. 

Now consider the additional constraints $x^0 \in D$ and $x^K \in E$, where $D \in \mathbb{R}^{n_0}$ and $E \in \mathbb{R}^{n_K}$ are polytopes. We then have that 
\begin{equation}
\{(x^0, x^K) : x^0 \in D, x^K \in E, (x, y, z) \in F(\Theta, B) \} \subseteq F(\Theta, B)
\end{equation}
The reason for making this, perhaps obvious point, is that a DNN programmed by \eqref{eq:nn-milp} may be part of a larger optimization problem with additional constraints on $x^0$ and $x^K$. These constraints may tighten variable bounds and should be considered when performing bound tightening.

The previous works on bound tightening (BT) of \eqref{eq:nn-milp}, listed in the introduction, did not include output bounds $E$. The primary goal has been to search in $F(\Theta, B)$ subject only to input bounds $D$, which explains the focus on \emph{forward-propagating} BT procedures. To find optimal bounds with the addition of output bounds $E$, it is required that the output bounds are propagated backwards through the DNN. In the following, we study several bound tightening procedures for \eqref{eq:nn-milp} with $D$ and $E$ included. We will highlight the procedures that are capable of propagating bounds backwards and subsequently investigate the effect of backwards propagation of bounds. We denote the directions of bound propagation by forward bound propagation (FBP) and backward bound propagation (BBP), respectively.

To simplify the notation, we assume in the following that $D$ and $E$ are hyperrectangles included in $B^0 = [L^{0}, U^{0}]$ and $B^K = [L^{K}, U^{K}]$.

\subsection{Relaxations of ReLU networks}
\label{sec:relu-relaxations}
We define the feasible set of the LP relaxation of \eqref{eq:nn-milp} as
\begin{equation}
R_R(\Theta, B) := \left\{ (x, s, z) : \eqref{eq:nn-milp-input-bounds}, \eqref{eq:nn-milp-linear}, \eqref{eq:nn-milp-bigm}-\eqref{eq:nn-milp-output-bounds} \text{ are satisfied}, 0 \leq z \leq 1 \right\},
\label{eq:feasible-set-lp-relaxation}
\end{equation}
where the binary restrictions on $z$ are replaced by $0 \leq z \leq 1$. This relaxation can be interpreted as a relaxation of the ReLU logic, and we refer to it as \emph{ReLU relaxation}. The ReLU relaxation makes it possible to increase the output from the ReLU and we get the following bounds on the output from node $j$ of layer $k$:
\begin{equation}
(W^{k} x^{k-1} + b^{k})_j \leq x_j^k \leq U_j^k \frac{(W^{k} x^{k-1} + b^{k})_j - L_j^k}{U_j^k - L_j^k}
\label{eq:relu-bound}
\end{equation}
We note two things from \eqref{eq:relu-bound}. First, even if $(W^{k} x^{k-1} + b^{k})_j < 0$ the unit can still give a positive output in the ReLU relaxation. Second, the upper bound depends on both $U_j^k$ and $L_j^k$, and a decrease in $U_j^k$ or increase in $L_j^k$ directly lowers the bound. This clearly motivates a study of BT procedures. 

We use the notation $R_R(\Theta, B, I)$ to denote a partial ReLU relaxation for which the ReLUs of nodes $(j,k)$ in the index set $I$ are relaxed. A MILP problem results from a partial ReLU relaxation. We let the LP relaxation $R_R(\Theta, B)$ correspond to $R_R(\Theta, B, I)$ with all nodes indexed by $I$.

When considering the bound tightening of a node $(j,k)$, it is possible to relax $P(\Theta, B)$ by removing constraints related to other layers than $k$ from \eqref{eq:nn-milp}. We may for example remove all layers following $k$ ($k+1, \ldots, K$) to simplify the problem while retaining FBP from layer $0$ through $k$. A weaker relaxation can be obtained by removing all layers but $k$ and $k-1$. This relaxation can be used to propagate bounds forward one layer at the time. In both these cases, where all layers following $k$ are excluded, node $(j,k)$ can be considered an output node. We may then also remove the constraints of all the other nodes in layer $k$ to reduce the problem size. We name this type of relaxation \emph{layer relaxation} and denote it by $R_L(\Theta, B, I)$, where we specify by the index set $I$ which nodes to remove/relax.

For the relaxations discussed above, we have that $P(\Theta, B) \subseteq R_R(\Theta, B, I)$ and $P(\Theta, B) \subseteq R_L(\Theta, B, I)$ for any index set $I$. Thus, a BT procedure utilizing these relaxation produces \emph{valid} bounds; i.e., the tightened bounds do not diminish the set of input-output pairs $F(\Theta, B)$.


\subsection{Prototype for bound tightening procedures}
The bound tightening procedures we discuss subsequently are special cases of the prototype in Algorithm \ref{alg:prototype}.

\begin{algorithm}[ht]
\caption{Prototype for bound tightening procedures}
\label{alg:prototype}
\begin{algorithmic}
\Require Parameters $\Theta$, initial bounds $B^-$ and scheme $S$
\State $B \gets B^-$ 
\State $k \gets 0$ 
\For{$k \leq K$} \Comment{Iterate over layers}
	\State $j \gets 1$ 
	\For{$j \leq n_k$} \Comment{Iterate over nodes}
		\State Build constraint set $C_{j}^{k}$ from $\Theta$ and $B$ following $S$
		\State Solve for upper bound $u_j^k = \max t_j^k$ s.t. $C_{j}^{k}$
		\State Solve for lower bound $l_j^k = \min t_j^k$ s.t. $C_{j}^{k}$
		\State Update bounds $B$: $U_j^k \gets u_j^k$ and $L_j^k \gets l_j^k$
		\State $j \gets j + 1$
	\EndFor
	\State $k \gets k + 1$
\EndFor \\
\Return Tightened bounds $B$
\end{algorithmic}
\end{algorithm}

Given a network parameterized by $\Theta$, some initial bounds $B^-$ and a scheme $S$, the prototypical procedure computes tightened bounds $B$. The scheme $S$ state for each node $(j,k)$ which constraint set $C_j^k$ that should be used. $C_j^k$ represents a subset of the constraints in \eqref{eq:nn-milp}, or a valid relaxation of these.

Bounds $l_{j}^{k}$ and $u_{j}^{k}$ of each node $(j,k)$ are computed so that $l_{j}^{k} \leq t_{j}^{k} \leq u_{j}^{k}$, by minimizing and maximizing the variable $t_{j}^{k}$ defined as:
\begin{equation}
t_{j}^{k} := \begin{cases}
x_{j}^{k} & k = 0 \\
x_{j}^{k} - s_{j}^{k} & k = 1, \ldots,K-1 \\
x_{j}^{k} & k = K \\
\end{cases}
\end{equation} 

If a strictly positive lower bound is found for a hidden node $(j,k)$, i.e. $L_j^k > 0$, the corresponding ReLU may be relaxed by fixing $z_j^k = 1$ and $s_j^k = 0$ (or removing them from the problem). Similarly, if $U_j^k < 0$, we may fix $z_j^k = 0$ and $x_j^k = 0$; for this case we say that the neuron $(j,k)$ is dead. This effectively reduces the number of binary variables in the formulation.

When tightening the bounds of \eqref{eq:nn-milp}, we should consider in which order to processes the nodes. For FBP, it is natural to process nodes in the order of layers, beginning with the nodes in the first hidden layer ($k=1$). However, when bound information can flow backwards from the output bounds $B^K$, it is not obvious in which order the nodes should be processed. In \ref{app:forward-vs-backward-propagation}, we provide an argument for using the same strategy also when bound information may propagate backwards.\footnote{These procedures must not be confused with the forward and backward propagation used in inference and training of neural networks.} The argument is that, since BBP is much weaker than FBP, the bound tightening procedure should focus on propagating bounds forward. In line with this argument, Algorithm \ref{alg:prototype} makes one forward pass through the layers, processing the nodes of each layer in the order they are indexed. We may apply the argument again when considering the ordering within each layer. Since the only interaction between the nodes in a layer is via BBP on subsequent nodes, different processing orders will likely perform similarly.

We note that $k$ is initialized to 0, since the input bounds $[L^0, U^0]$ may be tightened by the output bounds when the constraint sets $C_j^k$ allow for BBP. For instances of the procedure for which BBP cannot occur due to the choice of constraint set $C_j^k$, $k$ may be initialized to 1.

\subsection{Feasibility-based bound tightening (FBBT)}
FBBT subsumes BT procedures that use purely primal feasibility arguments and remove parts of the domains in which no feasible solutions are contained \cite{Gleixner2017}. These procedures rely on interval arithmetic to compute the bounds on constraint activations over the variable domains, and conversely, propagating the bounds on the constraint activities back to the variable domains.

Using simple interval arithmetic we may infer bounds $l_j^k \leq t_j^k \leq u_j^k$ for node $(j, k)$ in layer $k \geq 2$ as follows:
\begin{equation}
\begin{aligned}
u_j^k &= \sum\limits_{i=1}^{n_{k - 1}} \max \left\{ w_{ji}^{k} \max \{0, U_{i}^{k-1}\}, w_{ji}^{k} \max \{0, L_{i}^{k-1}\} \right\} + b_{j}^{k}, \\
l_j^k &= \sum\limits_{i=1}^{n_{k - 1}} \min \left\{ w_{ji}^{k} \max \{0, U_{i}^{k-1}\}, w_{ji}^{k} \max \{0, L_{i}^{k-1} \} \right\} + b_{j}^{k},
\end{aligned}
\label{eq:LRR-equations}
\end{equation}
where we have utilized the bounds on $x_j^{k-1}$ in \eqref{eq:bound-update-rule}. The same bounds can be found by solving the LP problems:
\begin{equation}
\begin{aligned}
u_j^k &= \max \left\{ t_{j}^{k} : t_{j}^{k} \in C_{j}^{k} \right\}, \\
l_j^k &= \min \left\{ t_{j}^{k} : t_{j}^{k} \in C_{j}^{k} \right\},
\end{aligned}
\label{eq:LRR-problem}
\end{equation}
for the constraint set 
\begin{equation}
C_{j}^{k} = \left\{ t_{j}^{k} : t_{j}^{k} = w_j^k x^{k-1} + b^k, x^{k-1} \in [\max\{0, L^{k-1}\}, \max\{0, U^{k-1}\}] \subset \mathbb{R}^{n_{k-1}} \right\}.
\label{eq:LRR-constraint-set}
\end{equation}
To compute the bounds in the first hidden layer $k=1$, we remove the inner max-operators in \eqref{eq:LRR-equations} or \eqref{eq:LRR-constraint-set}, since there is no ReLU operator on the input layer.

For a DNN, we may propagate the input bounds forward by solving \eqref{eq:LRR-equations} or \eqref{eq:LRR-problem}, and then update the bounds according to \eqref{eq:bound-update-rule} for units in successive layers, beginning with the units in layer $k = 1$. That is, we follow the prototypical procedure in Algorithm \ref{alg:prototype}, but initialize with $k=1$. This procedure employs both layer and ReLU relaxations, as seen from $C_j^k$, and we refer to it as LRR in the rest of this paper. The procedure is computationally cheap since it considers only the constraints in \eqref{eq:nn-milp-linear} for the node being tightened. 

In general, feasibility-based procedures compute bounds that are sub-optimal, due to the weak relaxations. Stronger BT procedures can be devised by including more constraints and exploiting integer information. For example, the LRR procedure propagates bounds in the forward direction only, and may not take advantage of the output bounds $B^K$. To do so, a BT procedure must be able to perform BBP, either by starting at the last layer, or by including enough constraints to link the variable domain being tightened to the output bounds. Next, we consider some alternative BT procedures that use stronger relaxations, but are more computationally demanding.

\subsection{Optimization-based bound tightening (OBBT)}
We devise several OBBT procedures for ReLU networks using the relaxations in Section \ref{sec:relu-relaxations}. These procedures rely on solving a series of $2\sum_{k=0}^{K}n_k$ optimization problems to find the tightest variable bounds on relaxations of \eqref{eq:nn-milp}.

\paragraph{RR procedure} First, we consider an OBBT procedure using the LP relaxation $R_R(\Theta, B)$ of \eqref{eq:nn-milp}. The procedure computes bounds on $t_j^k$ by solving the optimization problems:
\begin{align*}
u_j^k &= \max \left\{ t_j^k : (x, s, z) \in R_R(\Theta, B) \right\}, \\
l_j^k &= \min \left\{ t_j^k : (x, s, z) \in R_R(\Theta, B) \right\},
\end{align*}
for $k=0,\ldots,K$ and $j=1,\ldots,n_k$. Setting $C_j^k = R_R(\Theta, B)$ in Algorithm \ref{alg:prototype}, we obtain an LP-based OBBT procedure which we call RR. 


The subsequent procedures require solving $2\sum_{k=0}^{K}n_k$ MILP problems of increasing complexity (the number of constraints and binary variables increase with layer depth).

\paragraph{LR procedure} Next, we consider a procedure that utilizes layer relaxations $R_L(\Theta, B, I_j^k)$. When solving for node $(j,k)$, we construct the index set $I_j^k$ so that the constraints are removed for all nodes except $(j,k)$ and nodes in preceding layers $\{0,\ldots,k-1\}$. In the framework of Algorithm \ref{alg:prototype}, we initialize from $k=1$ and set $C_j^k = R_L(\Theta, B, I_j^k)$. We name this procedure LR since it only utilizes layer relaxations. 

\paragraph{SEMI-RR procedure} The LR procedure removes all layers following a node $(j,k)$ and may not propagate bounds backwards. To address this potential weakness, we devise a procedure that keeps all following layers, but relaxes their ReLU operators. The procedure is given by Algorithm \ref{alg:prototype} with $C_j^k = R_R(\Theta, B, I_j^k)$, where $I_j^k$ indexes the nodes in layers $\{k,\ldots,K-1\}$ (excluding layer $K$ since it is affine). The procedure is named SEMI-RR since it relaxes the ReLUs in layer $k$ onwards.

\paragraph{NO-R procedure}
Finally, we may compute the tightest possible bounds of a ReLU network by using the full constraint set in \eqref{eq:nn-milp}; i.e. we set $C_j^k = P(\Theta, B)$. The procedure does not relax any constraints and we therefore name it NO-R.

\subsection{Summary of bound tightening procedures}
The various BT procedures discussed above employ two different types of relaxations. They either relax the integrality constraint on the binary variables $z$ (ReLU relaxation), or remove nodes or layers and their respective constraints (layer relaxation), or utilize both types of relaxations. The procedures, which we have named based on the relaxations they employ, are listed in Table \ref{tab:bt-procedures}.

\begin{table}[ht]
\footnotesize
\caption{Bound tightening procedures.}
\begin{tabularx}{\textwidth}{lll*3{>{\raggedleft\arraybackslash}X}}
\toprule
BT procedure & BT type & Subproblem class & Layer relaxation & ReLU relaxation & Backward propagation** \\
\midrule 
LRR & FBBT & LP/interval arithmetic & \checkmark & \checkmark & -- \\
RR & OBBT & LP & -- & \checkmark & \checkmark \\
LR & OBBT & MILP & \checkmark & -  & -- \\
SEMI-RR & OBBT & MILP & -- & * & \checkmark \\
NO-R & OBBT & MILP & -- & -- & \checkmark \\
\bottomrule \noalign{\smallskip}
\multicolumn{5}{l}{*ReLUs of preceding nodes are not relaxed.} \\
\multicolumn{5}{l}{**BT procedure can propagate output bounds backwards through the network.}
\end{tabularx}
\label{tab:bt-procedures}
\end{table}

Let $B^{\text{LRR}}$ be the bounds found by running LRR, $B^{\text{RR}}$ by RR, and so on. By inspecting Table \ref{tab:bt-procedures} we see that $\text{MAD}(B^{\text{RR}}) \leq \text{MAD}(B^{\text{LRR}})$. Likewise, we have that $\text{MAD}(B^{\text{NO-R}}) \leq \text{MAD}(B^{\text{SEMI-RR}}) \leq \text{MAD}(B^{\text{LR}})$, given that subproblems are not time limited. In the absence of output bounds, we have that $\text{MAD}(B^{\text{NO-R}}) = \text{MAD}(B^{\text{SEMI-RR}}) = \text{MAD}(B^{\text{LR}})$, and LR is clearly the cheaper procedure. It is not trivial to compare the LP-based procedures with MILP-based procedures since they employ different relaxations. However, we do have that $\text{MAD}(B^{\text{SEMI-RR}}) \leq \text{MAD}(B^{\text{RR}})$, again assuming that subproblems are not time limited.

Several of the presented BT procedures occur in other works on verification of neural networks; namely, the LRR, RR, and LR procedures were employed in \cite{Tjeng2017,Cheng2017,Fischetti2018} for forward propagation of bounds. While NO-R simply uses the full constraint set in \eqref{eq:nn-milp}, the SEMI-RR has not to our knowledge been presented before.

We also note that multiple invocations of the RR or SEMI-RR procedure may strengthen the bounds. The additional tightening is possible due to the BBP ability of these procedures. I.e., the computation of a node's bounds is affected by the tightened bounds on deeper nodes from earlier runs. The effect was demonstrated in \cite{Tjeng2017} for the RR procedure.

\begin{remark}[Pre-computing bounds] Due to their computational cost, we only invoke the BT procedures once before optimization (or at the root node of the branch-and-bound tree). The bounds can then be stored and reused in subsequent optimizations. We rely on the solver's bound tightening capabilities further down the branch-and-bound tree.
\end{remark}

\begin{remark}[Subproblem time limit] For the MILP-based procedures LR, SEMI-RR, and NO-R we may limit the solution time of each subproblem to reduce the overall computational cost. For subproblems that do not terminate within this time limit we use the best bound found by the solver, ensuring that the computed bounds are still valid. We introduce a naming convention where we postfix the BT procedure name by the subproblem solution time limit (in seconds). For example, we write \emph{LR(60)} for the LR procedure with the solution time of the subproblems limited to 60 seconds.
\end{remark}

\begin{remark}[Bound initialization]
We initialize all OBBT procedures with bounds computed by LRR. This initialization is cheap and gives the procedures an identical starting point. Furthermore, it allows us to compare the procedures with the mean relative distance (MRD) statistic in \ref{app:bound-statistics}.
\end{remark}

\section{Numerical study}
\label{sec:numerical-study}
We provide a computational study to evaluate the feasibility of using ReLU networks as surrogate models in mixed-integer programs. We focus primarily on solution times for different network architectures and bound tightening procedures in our investigation. The BT procedures in Section \ref{sec:bound-tightening} are compared based on computational efficiency and the statistics in \ref{app:bound-statistics}, as well as optimization solution time. 

We first evaluate the BT procedures' ability to propagate output bounds backwards on a set of randomly generated neural networks. Next, we solve a series of increasingly challenging optimization problems to test the practical performance of the big-M formulation in \eqref{eq:nn-milp} with bounds computed by the proposed BT procedures. We finally solve an oil production optimization problem including ten ReLU network surrogate models. 

The neural networks are initialized and trained using TensorFlow \cite{tensorflow2015-whitepaper}. We solve all optimization problems using Gurobi 8.1 with default settings \cite{Gurobi}, on a machine equipped with an Intel Core i7-8700K processor and 32 GB of RAM memory.

\subsection{Bound tightening for randomly initialized ReLU networks with output bounds}
\label{sec:bound-tightening-with-output-bounds}
\begin{figure}[ht]
\centering
\subfloat[]{
\includegraphics[width=0.47\textwidth]{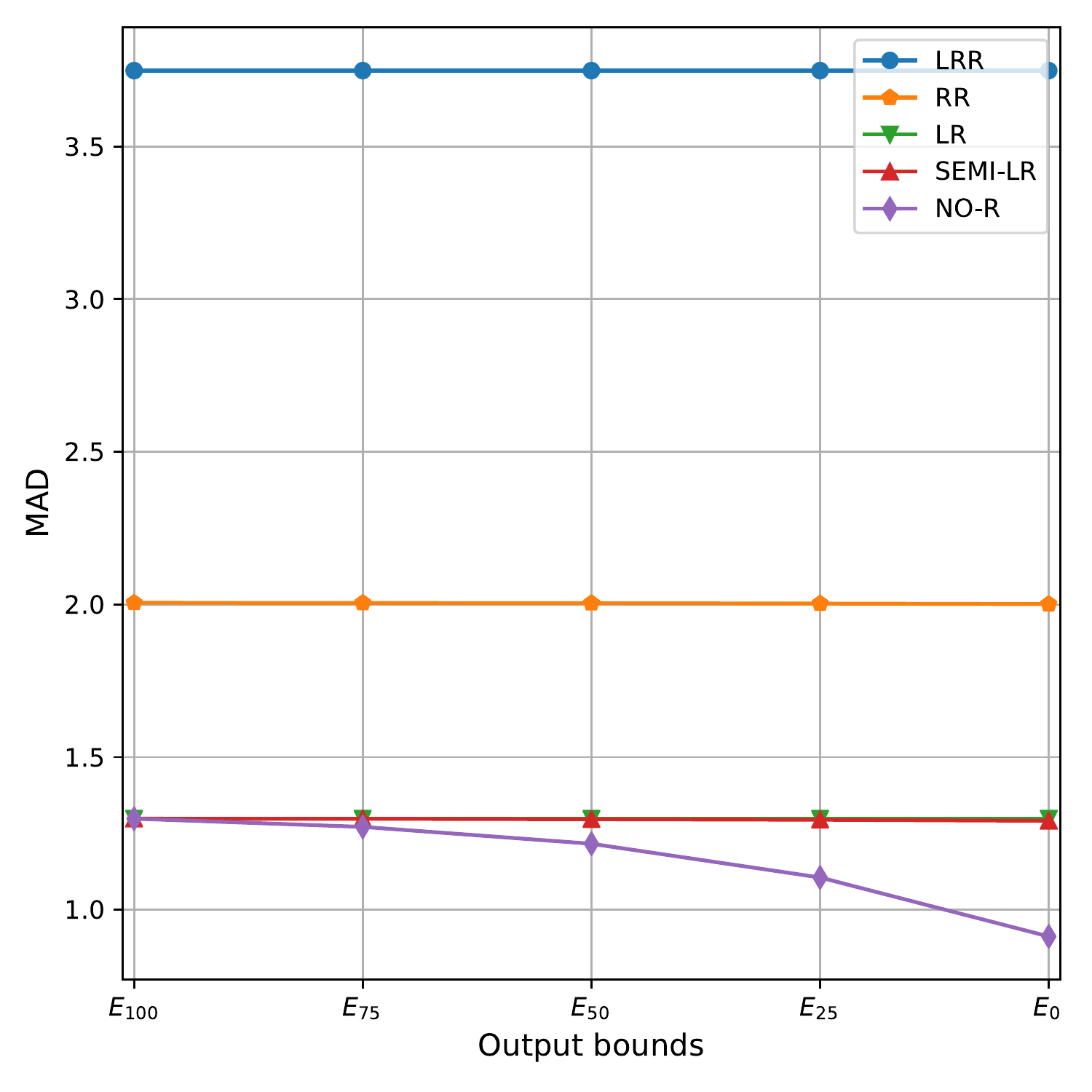}
\label{fig:bt_output_bounds_mad}
}
\hfill
\subfloat[]{
\includegraphics[width=0.47\textwidth]{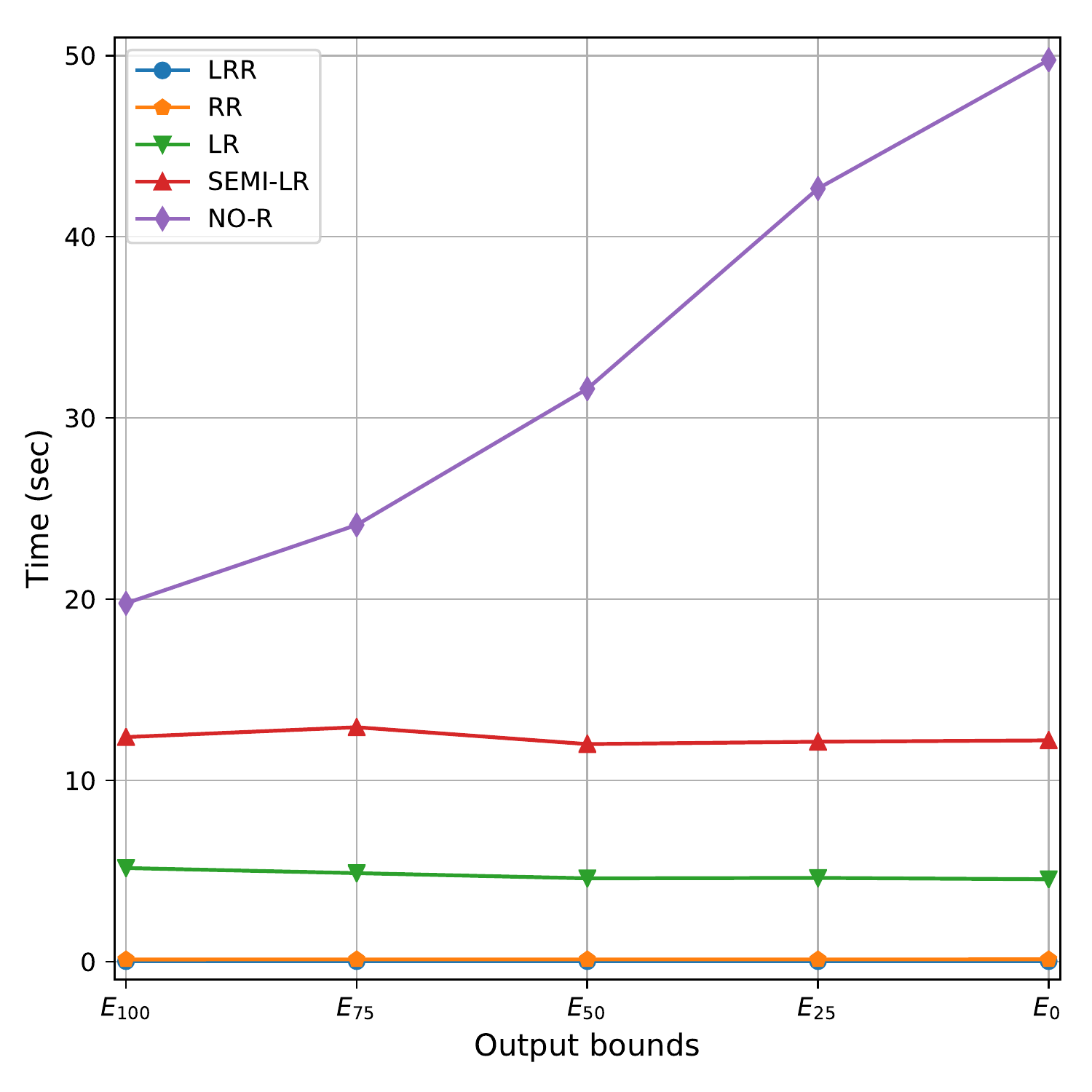}
\label{fig:bt_output_bounds_time}
}
\caption{Average MAD (a) and solution time (b) resulting from running the various BT procedures on ten neural networks. NO-R is the only procedure that significantly reduces the MAD as the output bounds are diminished.}
\label{fig:bt_output_bounds}
\end{figure}

We generate ten networks with layers $(3, 20, 20, 10, 1)$ to investigate the BT procedures' ability to propagate output bounds backwards. The networks parameters are initialized using the method by \cite{He2015}, so that for inputs $x^0 \sim \mathcal{N}(0, 1)$, the output $x^K \sim \mathcal{N}(0, 1)$. We set the input bounds to $B^0 = [-1, 1]^3$ and output bounds to $B^K = E_{100} := [-1, 1]$. We then proceed by running the procedures for diminishing output bounds: $E_{75} = [-0.75, 0.75]$, $E_{50} = [-0.5, 0.5]$, $E_{25} = [-0.25, 0.25]$, and $E_{0} = [0, 0] = \{0\}$. The results are shown in Figure \ref{fig:bt_output_bounds}, where time and MAD values are averaged over the ten generated networks.

Comparing the MAD values in Figure \ref{fig:bt_output_bounds_mad} for $E_{100}$, we see that LRR produces the loosest bounds, followed by RR. The MILP-based BT procedures, LR, SEMI-LR, and NO-R, obtain the same average MAD value. This is due to the fact that the output bounds $E_{100}$ have no tightening effect on preceding units; that is, $f(x) \in E_{100}$ for all $x \in B^0$. However, as the output bounds are diminished, we see that the average MAD decreases for NO-R. For comparison, we report the average MAD values for $E_0$ and $E_{100}$ (denoted $\text{MAD}_0$ and $\text{MAD}_{100}$) in Table \ref{tab:output-bounds-results}. As reported in the table, NO-R is the only BT procedure capable of significantly reducing the MAD by propagating the $E_{0}$ bounds backwards.

\begin{table}[ht]
\footnotesize
\caption{Average MAD for output bounds $E_0$ and $E_{100}$.}
\begin{tabularx}{\textwidth}{l*3{>{\raggedleft\arraybackslash}X}}
\toprule
BT procedure & $\text{MAD}_{100}$ & $\text{MAD}_{0}$ & $100 \times \text{MAD}_{0} / \text{MAD}_{100}$ \\
\midrule 
LRR & 3.65616 & 3.65616 & 100.00 \\
RR & 1.94224 & 1.91441 & 98.57 \\
LR & 1.25684 & 1.25684 & 100.00 \\
SEMI-LR & 1.25684 & 1.24049 & 98.70 \\
NO-R & 1.25684 & 0.95079 & \textbf{75.65} \\
\bottomrule
\end{tabularx}
\label{tab:output-bounds-results}
\end{table}

Turning to Figure \ref{fig:bt_output_bounds_time}, we see that the average solution times of the BT procedures are ordered as expected. The solution time increases with the number of constraints included in the BT procedure, and jumps considerably as binary variables are included. The ordering is the same as in the legend of the figure, LRR having the lowest and NO-R the highest solution times on average. The most interesting observation in this figure is that the solution times for NO-R increases as the output bound is tightened, while the other BT procedures seem to be unaffected. Possible explanations are that NO-R, to an increasing degree, must propagate the output bounds backwards, thus linking more variables (nodes) in the MILP problems and that finding feasible integer solutions are harder with a tight output bound. Furthermore, the constraints through which the output bounds must be propagated backwards are likely quite weak (see discussion in \ref{app:forward-vs-backward-propagation}), which may affect the numerical performance of the MILP solver.

\subsection{Optimization of approximated n-dimensional quadratic functions}
\label{sec:quadratic-approximations}
We consider $n$-dimensional quadratic functions 
\begin{align*}
q(x) = x^\top A x + b^\top x + c,
\end{align*}
where $x \in \mathbb{R}^n$, $A \in \mathbb{R}^{n \times n}$, $b \in \mathbb{R}^n$, and $c \in \mathbb{R}$. The coefficients are drawn as $A_{ij} \sim \mathcal{N}(0, 5)$, $b_{i} \sim \mathcal{N}(0, 1)$, and $c \sim \mathcal{N}(0, 1)$. The quadratic functions resulting from this construction are likely to be indefinite and non-separable, and thus difficult to optimize. Furthermore, the higher variance on the coefficients of the quadratic terms increases the curvature of the generated functions. Due to the curvature, many pieces are required to obtain an accurate piecewise-linear approximation of these functions.

\begin{figure}[ht!]
\begin{center}
\includegraphics[width=0.7\textwidth]{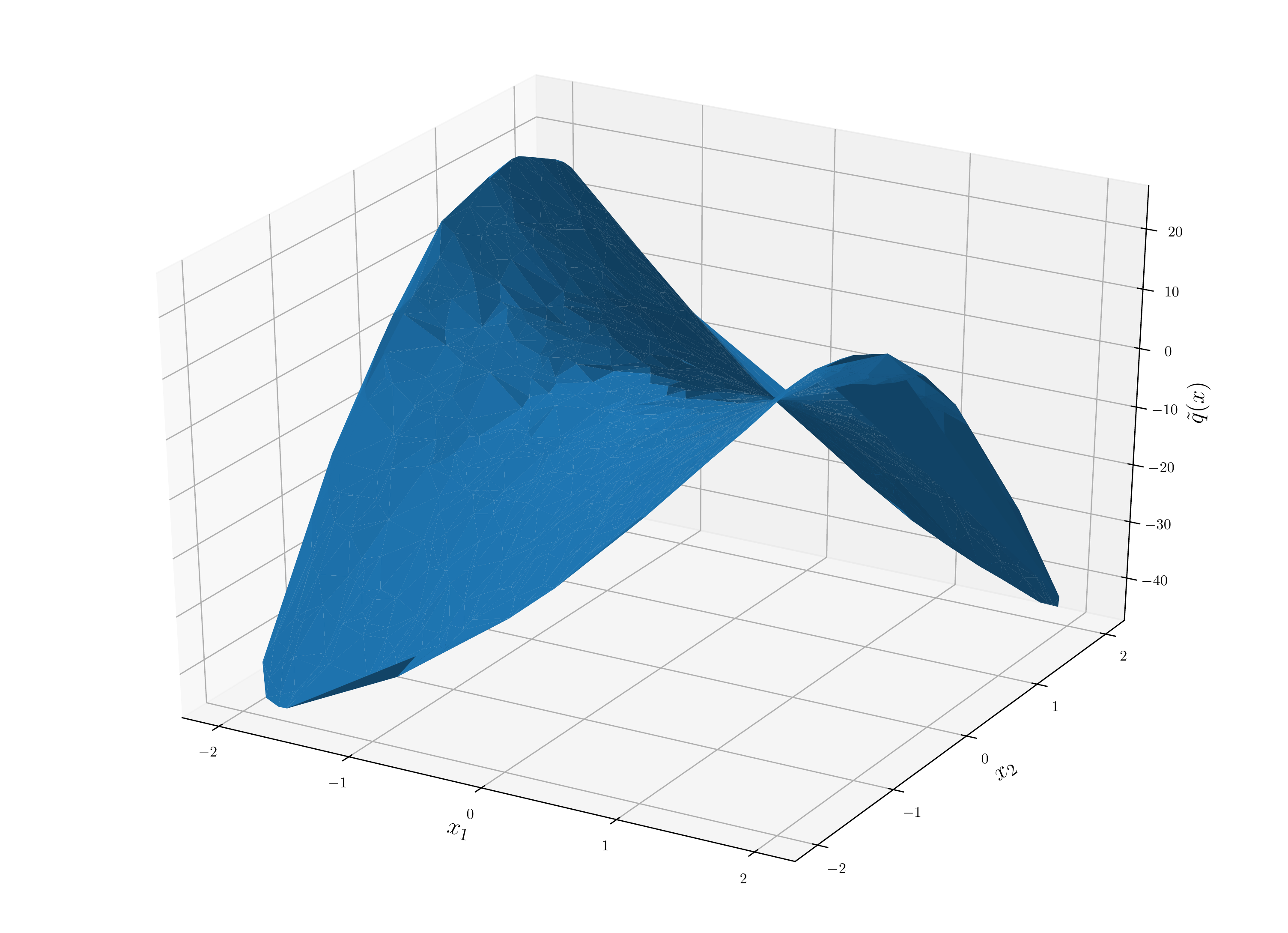}
\caption{ReLU network approximation of a quadratic function with an indefinite $A$ matrix.} 
\label{fig:quad-approx}
\end{center}
\end{figure}

For $n = 1, \ldots, 6$, we generate ten quadratic functions and approximate them by ReLU networks with the configurations given in Table \ref{tab:quadratic-approximations}. The table also lists the number of training samples, and the resulting \emph{mean absolute percentage error} (MAPE) averaged over the ten models. One of the two-dimensional approximations is shown in Figure \ref{fig:quad-approx}. For the displayed case, the drawn $A$ matrix is indefinite, leading to a saddle surface.

\begin{table}[ht!]
\footnotesize
\caption{ReLU network approximations of quadratic functions.}
\begin{tabularx}{\textwidth}{l*4{>{\raggedleft\arraybackslash}X}}
\toprule
$n$ & Layers & \# parameters & \# training samples & MAPE (\%)\\
\midrule 
1 & (1, 10, 5, 1) & 81 & 100 & 2.3 \\  
2 & (2, 20, 10, 1) & 281 & 500 & 3.5 \\
3 & (3, 40, 20, 1) & 1 001 & 2 000 & 3.7 \\
4 & (4, 50, 20, 1) & 1 291 & 5 000 & 3.7 \\
5 & (5, 50, 30, 30, 1) & 2 791 & 10 000 & 2.6 \\
6 & (6, 80, 40, 40, 1) & 5 481 & 20 000 & 2.0 \\
\bottomrule \noalign{\smallskip}
\end{tabularx}
\label{tab:quadratic-approximations}
\end{table}

Using the ReLU approximations, we consider optimization problems on the following form:
\begin{equation}
\min \left\{ f_{\Theta_1}(x) :  f_{\Theta_2}(x) = \alpha, x \in [-1, 1]^n \right\},
\label{eq:quad-problem}
\tag{$Q_n$}
\end{equation}
where $f_{\Theta_1}$ and $f_{\Theta_2}$ are ReLU networks, and $\alpha$ is a real constant. For $n = 1, \ldots, 6$, we construct five optimization problems by pairing the ReLU networks described above. 

We set the objective function to $f_{\Theta_1}$, and optimize under the constraint that the second network, $f_{\Theta_2}$, must intersect a plane at level $\alpha$ at the solution. Notice that the feasible set, specified by the constant $\alpha \in \mathbb{R}$, is disconnected when $f_{\Theta_2}$ has multiple crossings of the $\alpha$-plane. This can be seen for the saddle surface in Figure \ref{fig:quad-approx}, for which the intersection to an $\alpha$-plane is a hyperbola (unless the $\alpha$-plane intersects the saddle point).

The optimization results are given in Table \ref{tab:solve-times-quadratic}. When no feasible solution is found within the total time for one or more of the networks, we mark this with $\infty$ in the \emph{Gap} column. From the table, we see that the LRR procedure has the best performance for the low dimensional problems ($n = 1,2,3$). A likely explanation is that the solver's internal BT procedures work well for small network sizes, and that the stronger bound tightening procedures lead to computational overhead. For $n=4$, the RR procedure takes the lead by striking a good balance between time spent on BT and resulting bound tightness. For the largest problems ($n=5,6$), the LR(1) procedure is the best performer in terms of total solution time. The trend seems to be that larger network sizes benefit more from stronger bound tightening. We note that for $n=6$, we are able to solve all five instances when using the LR(1) procedure. Even with subproblem solution time is limited to one second, the SEMI-LR(1) and NO-R(1) procedures use several hundred seconds on bound tightening, leaving little time for solving the optimization problems within the total time limit of 600 seconds.

\begin{table}[ht!]
\footnotesize
\caption{Average solution times and optimality gap for \ref{eq:quad-problem} (five instances for each value of $n$). Solution times for bound tightening ($T_{\text{BT}}$) and optimization ($T_{\text{OPT}}$) are given in seconds. Total solution time, $T_{\text{BT}} + T_{\text{OPT}}$, is limited to 600 seconds.}
\begin{tabularx}{\textwidth}{ll*5{>{\raggedleft\arraybackslash}X}}
\toprule
$n$ & BT procedure & $T_{\text{BT}}$ & $T_{\text{OPT}}$ & $T_{\text{BT}} + T_{\text{OPT}}$ & Gap (\%) & \# solved \\
\midrule 
1 & LRR & 0.013 & 0.003 & \textbf{0.016} & 0 & 5/5 \\
  & RR & 0.026 & 0.002 & 0.028 & 0 & 5/5 \\
  & LR(1) & 0.054 & 0.003 & 0.057 & 0 & 5/5 \\
  & SEMI-LR(1) & 0.138 & 0.002 & 0.140 & 0 & 5/5 \\
  & NO-R(1) & 0.070 & 0.001 & 0.071 & 0 & 5/5 \\
  \noalign{\smallskip\smallskip}
2 & LRR & 0.024 & 0.070 & \textbf{0.094} & 0 & 5/5 \\
  & RR & 0.078 & 0.054 & 0.132 & 0 & 5/5 \\
  & LR(1) & 0.374 & 0.042 & 0.416 & 0 & 5/5 \\
  & SEMI-LR(1) & 1.575 & 0.047 & 1.622 & 0 & 5/5 \\
  & NO-R(1) & 1.959 & 0.037 & 1.996 & 0 & 5/5 \\
  \noalign{\smallskip\smallskip}
3 & LRR & 0.065 & 0.418 & \textbf{0.48} & 0 & 5/5 \\
  & RR & 0.293 & 0.302 & 0.60 & 0 & 5/5 \\
  & LR(1) & 3.772 & 0.237 & 4.01 & 0 & 5/5 \\
  & SEMI-LR(1) & 17.42 & 0.235 & 17.66 & 0 & 5/5 \\
  & NO-R(1) & 36.65 & 0.206 & 36.86 & 0 & 5/5 \\
  \noalign{\smallskip\smallskip}  
4 & LRR & 0.082 & 3.497 & 3.58 & 0 & 5/5 \\
  & RR & 0.403 & 1.644 & \textbf{2.05} & 0 & 5/5 \\
  & LR(1) & 7.717 & 1.217 & 8.93 & 0 & 5/5 \\
  & SEMI-LR(1) & 32.02 & 1.273 & 33.29 & 0 & 5/5 \\
  & NO-R(1) & 84.29 & 0.528 & 84.82 & 0 & 5/5 \\
  \noalign{\smallskip\smallskip}  
5 & LRR & 0.164 & 204.1 & 204.3 & 10.6 & 4/5 \\
  & RR & 2.252 & 137.9 & 140.1 & 0.4 & 4/5 \\
  & LR(1) & 95.59 & 13.2 & \textbf{108.8} & 0 & 5/5 \\
  & SEMI-LR(1) & 181.2 & 17.7 & 198.9 & 0 & 5/5 \\
  & NO-R(1) & 303.9 & 16.0 & 319.9 & 0 & 5/5 \\
  \noalign{\smallskip\smallskip}
6 & LRR & 0.286 & 433.0 & 433.3 & $\infty$ & 2/5 \\
  & RR & 5.97 & 295.0 & 301.0 & $\infty$ & 3/5 \\
  & LR(1) & 195.8 & 78.1 & \textbf{273.9} & 0 & 5/5 \\
  & SEMI-LR(1) & 335.0 & 128.2 & 463.2 & $\infty$ & 3/5 \\
  & NO-R(1) & 472.9 & 69.0 & 541.9 & 2.7 & 3/5 \\
\bottomrule \noalign{\smallskip}
\end{tabularx}
\label{tab:solve-times-quadratic}
\end{table}

The results given in Table \ref{tab:bt-results-quadratic}, \ref{app:bound-tightening-results}, show that the BT procedures that utilize stronger relaxations also compute tighter bounds, as expected. Another expected observation is that the number of dead neurons identified is strongly correlated with bound tightness.

In terms of BT solution times, we see that the interval arithmetic-based LRR procedure scales well with network size. The LP-based procedures scale moderately, while the MILP-based procedures scale poorly and quickly become computationally demanding. For the larger networks, the subproblem time limit of one second helps to limit BT solution times, at the cost of looser bounds. This explains why LR, which has the smallest subproblems among the MILP-based procedures, computes the tightest bounds for $n=6$.

\subsection{Oil production optimization case} 
\label{sec:prodopt-case}
To test the practical performance of the proposed methods we solve a production optimization case from an offshore petroleum production field involving eight subsea wells producing oil, gas, and water. Each of the wells produce to a topside processing facility, a platform, via one of two risers (pipelines transporting the fluid from the seabed to the topside processing facility). 

An offshore production system consists of several interconnected modules with corresponding interdependencies. The reservoir is the subsurface structure where oil and gas are located prior to extraction. A number of wells are drilled into the reservoir. The flow coming from a well is a mixture of water, oil and gas, called production phases, and is routed to a platform through a network of pipelines. Manifolds are used to connect different pipelines and to mix the incoming flows. At the platform, separators split the production phases and route the resulting streams of oil and gas to an export line. The export lines lead production phases off-site. The production flow from a well may be routed to a subset of the separators on the platform. A well can only be routed to a single separator at any given time. The routing decision must take into consideration the interdependencies between wells which produce to the same separator.

\begin{figure}[ht]
	\begin{center}
		\begin{tikzpicture}[
		scale=0.9,
		node/.style={circle,inner sep=0pt,minimum size=8mm,fill=gray!30,draw=black}]
		\node[node] (n1) at (0,0){1};
		\node[anchor=east] at (-0.4,0){$f_1$};
		\node[node] (n2) at (1.5,0){2};
		\node[node] (n3) at (3,0){3};
		\node[node] (n4) at (4.5,0){4};
		\node[node] (n5) at (6,0){5};
		\node[node] (n6) at (7.5,0){6};
		\node[node] (n7) at (9,0){7};
		\node[node] (n8) at (10.5,0){8};
		\node[node] (n9) at (4.5,2.5){9};
		\node[node] (n10) at (6,2.5){10};
		\node[node] (n11) at (4.5,5){11};
		\node[node] (n12) at (6,5){12};
		\draw[dashed, ->, black] (n1) to (n9);
		\draw[dashed, ->, black] (n1) to (n10);
		\draw[dashed, ->, black] (n2) to (n9);
		\draw[dashed, ->, black] (n2) to (n10);
		\draw[dashed, ->, black] (n3) to (n9);
		\draw[dashed, ->, black] (n3) to (n10);
		\draw[dashed, ->, black] (n4) to (n9);
		\draw[dashed, ->, black] (n4) to (n10);
		\draw[dashed, ->, black] (n5) to (n9);
		\draw[dashed, ->, black] (n5) to (n10);
		\draw[dashed, ->, black] (n6) to (n9);
		\draw[dashed, ->, black] (n6) to (n10);
		\draw[dashed, ->, black] (n7) to (n9);
		\draw[dashed, ->, black] (n7) to (n10);
		\draw[dashed, ->, black] (n8) to (n9);
		\draw[dashed, ->, black] (n8) to (n10);
		\draw[->, black] (n9) [anchor=east] to node {$g_{(9,11)}$} (n11);
		\draw[->, black] (n10) to (n12);
		\node[anchor=west] at (12,0){Wells, $\wellnodes$};
		\node[anchor=west] at (12,1.25){Pipelines, $\dedges$};
		\node[anchor=west] at (12,2.5){Manifolds, $\mannodes$};
		\node[anchor=west] at (12,3.75){Risers, $\risedges$};
		\node[anchor=west] at (12,5){Separators, $\sepnodes$};
		\end{tikzpicture}
		\caption{Production system flow graph, with nodes represented by grey circles and edges by arrows. Discrete edges are dashed. The nonlinearities $f_1$ and $g_{(9,11)}$, related to Node 1 and Edge (9,11), are shown.}
		\label{fig:flow-graph}
	
	\end{center}
\end{figure}
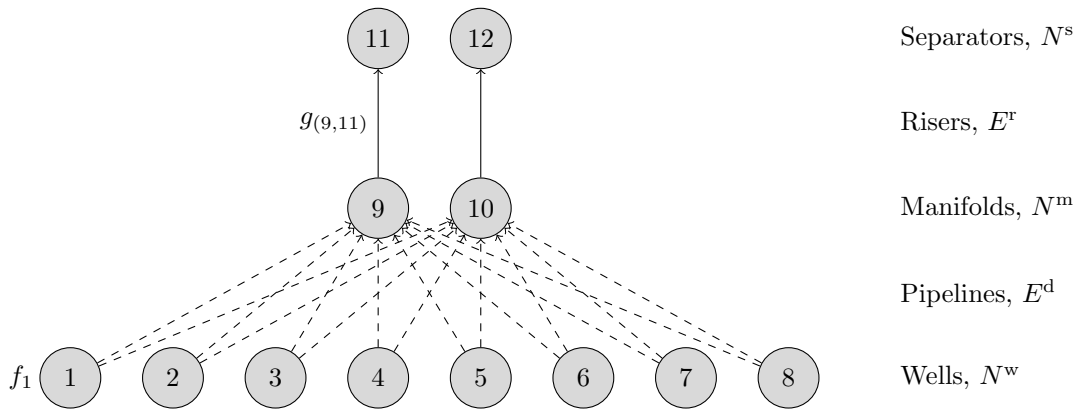

The production system is modeled by a directed acyclic graph $G=(\nodes, \edges)$, with nodes $\nodes$ and edges $\edges$. An illustration of the graph is given in Figure \ref{fig:flow-graph}. To simplify modeling, we use the utility sets in Table \ref{tab:utility-sets}. The set of nodes is $\nodes = \wellnodes \cup \mannodes \cup \sepnodes$, and the set of edges is $\edges = \dedges \cup \risedges$. As illustrated in Figure \ref{fig:flow-graph}, each well node ($i \in \wellnodes$) has two leaving \emph{discrete} edges. By allowing zero or one of these edges to be active, we model routing and on/off switching of well flows.

\begin{table}[ht]
\footnotesize
\caption{Utility sets}
\begin{tabularx}{\textwidth}{lX}
\toprule
Set & Description \\
\midrule 
$\nodes$ & Set of nodes in the network. \\
$\wellnodes$ & Set of \emph{well} (source) nodes in the network. $\wellnodes \subset \nodes$. \\
$\mannodes$ & Set of \emph{manifold} nodes in the network. $\mannodes \subset \nodes$. \\
$\sepnodes$ & Set of \emph{separator} (sink) nodes in the network. $\sepnodes \subset \nodes$. \\
$\edges$ & Set of edges in the network. An edge $e=(i,j)$ connects node $i$ to node $j$, where $i,j \in \nodes$. \\
$\dedges$ & Set of \emph{discrete} edges that can be \emph{open} or \emph{closed}. $\dedges \subset \edges$. \\
$\risedges$ & Set of \emph{riser} edges. $\risedges \subset \edges$. \\
$\edges_{i}^{\text{in}}$ & Set of edges entering node $i$, i.e. $\edges_{i}^{\text{in}} = \{ e : e=(j,i) \in \edges \}$. \\
$\edges_{i}^{\text{out}}$ & Set of edges leaving node $i$, i.e. $\edges_{i}^{\text{out}} = \{ e : e=(i,j) \in \edges \}$. \\
$\phases$ & $\phases = \{\oil,\gas,\water\}$, denoting the flow rate of oil, gas, and water, respectively. \\
\bottomrule
\end{tabularx}
\label{tab:utility-sets}
\end{table}

The variables are the pressure $p_i$ at each node $i \in \nodes$, the flow rate $q_{e,c}$ of phase $c \in \phases$ on edge $e \in \edges$, and the binary variable $y_e$ on the discrete edge $e \in \dedges$. The variables $y_e$ are used to model routing and on/off switching of wells. The problem also includes constants $c_{e,\text{gor}}$ (gas-oil ratio), $c_{e,\text{wor}}$ (water-oil ratio), flow rate bounds $q_{e,c}^{L}$ and $q_{e,c}^{U}$, pressure bounds $p_{i}^{L}$ and $p_{i}^{U}$, and separator pressures $p_{i}^{s}$, which are set to realistic values. For further details about the modeling approach we refer the reader to \cite{Grimstad2016b}. 

The complete formulation of the production optimization problem is given below.

\begin{subequations}
\begin{align}
	\underset{y,q,p}{\text{maximize}} \quad & z = \sum\limits_{e \in \risedges} q_{e,\oil} \label{eq:production-optimization-obj} \\ 
	\text{subject to}\quad & \sum\limits_{e \in \edges_{i}^{\text{in}}} q_{e,c} = \sum\limits_{e \in \edges_{i}^{\text{out}}} q_{e,c}, & \forall c \in \phases, i \in \mannodes \label{eq:production-optimization-balance} \\	
	& p_j = g_e(q_{e, \oil}, q_{e, \gas}, q_{e, \water}, p_i), & \forall e \in \risedges \label{eq:production-optimization-riser} \\	
	& (-p_j^U+p_i^L) (1-y_e) \leq p_i - p_j \leq (p_i^U - p_j^L) (1-y_e), & \forall e \in \dedges \label{eq:production-optimization-zero-drop} \\
	& \sum\limits_{e \in \edges_{i}^{\text{out}}} y_e \leq 1, & \forall i \in \wellnodes \label{eq:production-optimization-route} \\
	& y_{e}q_{e,c}^{L} \leq q_{e,c} \leq y_{e}q_{e,c}^{U}, & \forall c \in \phases, e \in \dedges \label{eq:production-optimization-quantity} \\
	& p_{i}^{L}	\leq p_{i} \leq p_{i}^{U}, & \forall i \in \nodes \label{eq:production-optimization-pressure}  \\
	& \sum\limits_{e \in \edges_{i}^{\text{out}}} q_{e, \oil} = f_{i}(p_i), &\forall i \in \wellnodes \label{eq:production-optimization-well} \\
	& \sum\limits_{e \in \edges_{i}^{\text{out}}} q_{e, \gas} = c_{e, \gor} \sum\limits_{e \in \edges_{i}^{\text{out}}} q_{e, \oil}, &\forall i \in \wellnodes  \label{eq:production-optimization-gor} \\
	& \sum\limits_{e \in \edges_{i}^{\text{out}}} q_{e, \water} = c_{e, \wor} \sum\limits_{e \in \edges_{i}^{\text{out}}} q_{e, \oil}, &\forall i \in \wellnodes  \label{eq:production-optimization-wor} \\
	& p_i = p_{i}^{s}, & \forall i \in \sepnodes \label{eq:production-optimization-const} \\	
	& y_e \in \{0,1\}, & \forall e \in \dedges  \label{eq:production-optimization-var}
\end{align}
\label{eq:production-optimization}
\end{subequations}

The objective function \eqref{eq:production-optimization-obj} maximizes the total oil production. Constraints \eqref{eq:production-optimization-balance} are mass balance constraints that make sure that the flow into a manifolds equals the flow out of it. The riser pressure drop functions $g_e(q_{e, \oil}, q_{e, \gas}, q_{e, \water}, p_i)$ for $e \in \risedges$ in equation \eqref{eq:production-optimization-riser} calculate the pressure at the separator as a function of the flow rates $(q_{e, \oil}, q_{e, \gas}, q_{e, \water})$ in the riser and the pressure at the start of the riser $p_i$. We assume no pressure drop over the discrete edges $(e \in \dedges)$ that are open ($y_e = 1$), as specified by the big-M constraints in \eqref{eq:production-optimization-zero-drop}. Constraints \eqref{eq:production-optimization-route} state that the flow from a well can be routed through at most one pipeline, while constraints \eqref{eq:production-optimization-quantity} connect flow rate through a pipeline with the routing decision. The flow must be within upper and lower bounds if the pipeline is open, and zero otherwise. The bounds on the pressure at each node are given by constraints \eqref{eq:production-optimization-pressure}. The well performance curves $f_i(p_i)$ for $i \in \wellnodes$ in equation \eqref{eq:production-optimization-well} relate the pressure at well $i$, $p_i$, with the flow rate of oil out of the well $q_{e,oil}$ for $e \in \edges_{i}^{\text{out}}$. Constraints \eqref{eq:production-optimization-gor} and \eqref{eq:production-optimization-wor} connect the flow rates of the different phases through the gas-oil and water-oil ratios. The pressure at the separators is fixed through constraints \eqref{eq:production-optimization-const} and constraints \eqref{eq:production-optimization-var} define the binary variables.

The optimization problem in \eqref{eq:production-optimization} contains a total of ten nonlinearities that we model with ReLU networks: eight well performance curves $f_i(p_i)$ for $i \in \wellnodes$, and two riser pressure drop functions $g_e(q_{e, \oil}, q_{e, \gas}, q_{e, \water}, p_i)$ for $e \in \risedges$. To test the effect of network depth on solution time, we consider the shallow and deep network architectures given in Table \ref{tab:prodopt-nn}. The networks are trained on scattered data sampled from a multiphase flow simulator. Each well and riser is modeled and simulated individually in the simulator. We simulate 50 data points for each well, and 4000 for each riser. Since the simulated data is free of noise, the DNNs are trained with a low $L_2$ penalty on the parameters. The ReLU network surrogates achieve a mean absolute percentage error (MAPE) of less than 1\%, meaning that they approximate the simulator with high accuracy.

\begin{table}[ht]
\footnotesize
\caption{Deep neural networks for wells and risers.}
\begin{tabularx}{\textwidth}{ll*3{>{\raggedleft\arraybackslash}X}}
\toprule
DNN & Layers & \# parameters & \# 0--1 var. & MAPE (\%) \\
\midrule 
Shallow well nets & (1, 20, 20, 1) & 481 & 40 & 0.73 \\
Shallow riser nets & (4, 50, 50, 1) & 2851 & 100 & 0.30 \\
\noalign{\smallskip\smallskip} 
Deep well nets & (1, 10, 10, 10, 10, 1) & 361 & 40 & 0.48 \\
Deep riser nets & (4, 20, 20, 20, 20, 20, 1) & 1801 & 100 & 0.42 \\
\bottomrule \noalign{\smallskip}
\end{tabularx}
\label{tab:prodopt-nn}
\end{table}


The MILP formulations of the shallow and deep networks in Table \ref{tab:prodopt-nn} have the same number of continuous and binary variables. There are 40 binary variables for the wells and 100 binary variables for the risers. The resulting optimization problem \eqref{eq:production-optimization} has a total of 536 binary variables (including 16 binary variables for well routing) and 1146 continuous variables. The number of constraints is 1416 with shallow networks, versus 1338 with deep networks. The problem size is thus similar with shallow and deep networks.

The BT procedures achieve the values reported in Tables \ref{tab:bt-results-shallow-well-nets}-\ref{tab:bt-results-deep-well-nets} for the well networks and Tables \ref{tab:bt-results-shallow-riser-nets}-\ref{tab:bt-results-deep-riser-nets} for the riser networks (see \ref{app:bound-tightening-results}). Due to the constraints on the separator pressures, the riser networks are output bounded to the singleton $\{p_i^s\}$ for $i \in \sepnodes$. Looking at the MRD values for the riser networks, it is evident that NO-R computes tighter bounds than the other BT procedures which are unable to fully utilize the output bounds via BBP. As can be seen from the solution times, the NO-R procedure is computationally expensive. For the shallow riser networks the average run time is 306 seconds, while for the deep riser networks it is 5080 seconds. Limiting the subproblem solution time to 60 seconds helps for the deeper networks, bringing the average solution time down to 1771 seconds for NO-R(60). The limitation in solution time do not seem to affect the bound tightness by much for these networks sizes, and NO-R(60) achieves a MRD of 0.1 \% compared to NO-R. For the well networks and the shallow riser networks, the subproblem time limit is never effective, which means that NO-R and NO-R(60) compute identical bounds.

With the tightened bounds, we solve the production optimization problem in \eqref{eq:production-optimization}. The results are given in Tables \ref{tab:shallow-prodopt-results} and \ref{tab:deep-prodopt-results} for the shallow and deep ReLU networks, respectively. The effect of bound tightness on optimization solution time is striking. With the shallow networks, we are only able to close the optimality gap within one hour using the bounds computed using the LR(60), NO-R(60) and NO-R procedure.  For the deep networks, the gap is closed when we use bounds computed by the NO-R(60) and NO-R procedure. Notice that the comparison is somewhat unfair since the overall computation time used by the NO-R procedure far exceeds that of the other BT procedures. It is, however, interesting to compare the optimization solution times ($T_{\text{OPT}}$) for the various BT procedures. 

With NO-R(60), we solve both the deep and shallow problem within one hour, which is acceptable for a practical application of \eqref{eq:production-optimization}. For larger networks with depth $K \geq 5$ and size exceeding 100 nodes, we are not able to solve \eqref{eq:production-optimization} within the practical time limit of one hour. 

With the optimal bounds computed by NO-R for the deep networks, the optimization converges in just 5 seconds. The expensive BT procedures thus show some merit on these problems, especially for applications where it is interesting to solve the problem many times; perhaps with small adjustments to the problem for each optimization run (required that the adjustments do not invalidate the computed bounds).

\begin{table}[ht]
\footnotesize
\caption{Solution times for production optimization problem with shallow ReLU networks. Bound tightening ($T_{\text{BT}}$) and optimization ($T_{\text{OPT}}$) solution times are given in seconds. $T_{\text{OPT}}$ is limited to 3600 seconds. $z^\star$ is the objective value of the best found solution and Gap is the difference between the upper and lower bounds relative to the lower bound.}
\begin{tabularx}{\textwidth}{l*5{>{\raggedleft\arraybackslash}X}}
\toprule
BT procedure & $T_{\text{BT}}$ & $T_{\text{OPT}}$ & $T_{\text{BT}} + T_{\text{OPT}}$ & $z^\star$ & Gap (\%) \\
\midrule 
LRR & 0.3 & 3600 & 3600 & 1.2864 & 7.5 \\
RR & 1.6 & 3600 & 3602 & 1.2864 & 4.7 \\
LR(60) & 18.6 & 2383 & 2402 & 1.2864 & 0 \\
SEMI-RR(60) & 115.0 & 3600 & 3715 & 1.2864 & 15.5 \\
NO-R(60) & 613.2 & 20 & \textbf{633} & 1.2864 & 0 \\
NO-R & 613.2 & 20 & \textbf{633} & 1.2864 & 0 \\
\bottomrule \noalign{\smallskip}
\end{tabularx}
\label{tab:shallow-prodopt-results}
\end{table}

\begin{table}[ht]
\footnotesize
\caption{Solution times for production optimization problem with deep ReLU networks. Bound tightening ($T_{\text{BT}}$) and optimization ($T_{\text{OPT}}$) solution times are given in seconds. $T_{\text{OPT}}$ is limited to 3600 seconds. $z^\star$ is the objective value of the best found solution and Gap is the difference between the upper and lower bounds relative to the lower bound.}
\begin{tabularx}{\textwidth}{l*5{>{\raggedleft\arraybackslash}X}}
\toprule
BT procedure & $T_{\text{BT}}$ & $T_{\text{OPT}}$ & $T_{\text{BT}} + T_{\text{OPT}}$ & $z^\star$ & Gap (\%) \\
\midrule 
LRR & 0.3 & 3600 & 3600 & -- & -- \\
RR & 2.5 & 3600 & 3603 & -- & -- \\
LR(60) & 220.1 & 3600 & 3820 & -- & -- \\
SEMI-RR(60) & 273.8 & 3600 & 3874 & -- & -- \\
NO-R(60) & 3543.2 & 57 & \textbf{3600} & 1.3049 & 0 \\
NO-R & 10161.2 & 5 & 10166 & 1.3049 & 0 \\
\bottomrule \noalign{\smallskip}
\end{tabularx}
\label{tab:deep-prodopt-results}
\end{table}

\begin{remark}
As an alternative approach to ReLU networks, we may model the well and riser nonlinearities using a traditional MILP formulation that interpolates the sample points, as discussed in Section \ref{sec:pwl-approximation}. For each nonlinearity, we use Delaunay triangulation on the samples to partition the domain into simplices. Specifically, for each riser model we create a Delaunay triangulation of 4000 scattered, four-dimensional data points. With the partitions, we obtain PWL approximations using the DLog model in \cite{Vielma2010}. The resulting formulation of \ref{eq:production-optimization} has 90 binary variables, 718 656 continuous variables, and 287 constraints. When attempting to solve this problem, the solver quickly runs out of memory without finding any feasible solution. We thus deem this approach unsuccessful for the production optimization case with irregularly sampled data points.
\end{remark}

\section{Concluding remarks}
\label{sec:concluding-remarks}


As we discussed in Section \ref{sec:pwl-approximation}, ReLU networks offer practitioners a versatile framework for modeling PWL functions. ReLU networks scale well to high input dimensions and can be trained on large datasets of scattered and noisy data. The MILP formulation in \eqref{eq:nn-milp} enables the embedding of ReLU networks in optimization problems that can be solved by state-of-the-art MILP solvers.

In statistical learning it is common to follow the principle of \emph{Occam's razor} and avoid unnecessary model complexity when searching for models that generalize well. According to this principle, deep ReLU networks are favorable to shallow networks since they can achieve the same complexity (number of linear partition regions) with fewer parameters. An example of this is given in Table \ref{tab:prodopt-nn}, where shallow (2851 parameters) and deep (1801 parameters) riser networks achieve similar accuracy. Smaller networks also tend to lower solution times in the optimization, which would lead us to conclude that deep networks are favorable for surrogate modeling. Paradoxically, the deepening of networks which has yielded favorable results on modeling tasks, seems to make the resulting MILP models harder to optimize (as shown in the numerical results). One explanation to this may be that with deeper networks, more variables are interlinked via the block-angular structure in Figure \ref{fig:ReLU-block-angular}. Thus, what is a strength in the modeling setting, may become a weakness in the optimization setting. We are thus compelled to search for a good trade-off between model depth and optimization efficiency. The same observation was reported in \cite{Schweidtmann2018,Fischetti2018} for global optimization of DNNs, supporting the claim that deeper networks are more challenging to optimize. To conclude, we observe, as others have before us and as theory predicts, that the feasibility of using the MILP formulation quickly fades with increasing network sizes. Our numerical study indicates that the MILP formulation is unfit for applications involving large ReLU networks with thousands of hidden nodes. For the oil production optimization case, we reached a practical limit on 100 hidden nodes for the deep network architectures. The bulk of the solution time is then used on bound tightening the large riser networks. 

Solution times are sensitive to tightness of the variable bounds in the MILP formulation, which directly translates to the big-M values in \eqref{eq:nn-milp}. Our study of bound tightening procedures, ranging from the cheap LRR procedure to the optimal NO-R procedure, shows that it is problem specific which procedure strikes the best computational efficiency. For optimization problems with small networks embedded, the cheaper procedures seem to perform best, while for larger networks the more expensive procedures perform best. Somewhat surprisingly, the MILP-based procedures, usually thought to be too computationally expensive, perform quite well on the more challenging problems. This is further amplified in applications where the optimization problem is solved many times using the same bounds. Here, spending the extra time in the more expensive bound tightening procedures can really pay off.

When a ReLU network is embedded in an optimization problem it is likely to be subject to output bounds. Our study is the first to investigate bound tightening of \eqref{eq:nn-milp} in the presence of output bounds. In Table \ref{tab:bt-procedures} we summarize the procedures and identify the ones capable of BBP. The numerical results in Section \ref{sec:numerical-study} show that ReLU relaxation significantly reduces bound tightness for both FBP and BBP. From the results, and the argument in \ref{app:forward-vs-backward-propagation}, we conclude that BBP is less effective than FBP. Of the studied procedures, only NO-R is capable of exploiting output bounds to reduce bounds; see Section \ref{sec:bound-tightening-with-output-bounds}. This may explain why NO-R is the best-performing procedure for the challenging oil production optimization case.

We finally remark that, while the application of \eqref{eq:nn-milp} is practically limited to small-sized networks, many interesting nonlinearities can be approximated by networks of this size (several examples can be found in works cited in the introduction). Our case studies show that the MILP formulation of ReLU networks in \eqref{eq:nn-milp} is an attractive approach to surrogate modeling in process optimization.

\subsection{Promising research directions}
Optimization of ReLU networks is currently an active area of research. Below, we highlight some promising research topics that may allow practitioners to embed larger ReLU networks in optimization problems.
\begin{itemize}
\item In this work we combined LRR with the other presented BT procedures. Other combinations may be more efficient in terms of tightening per processing time unit. For example, by combining RR with NO-R, we may first solve a round of LP subproblems, before we invoke the more expensive NO-R for further tightening. 
\item The feasible region defined by the set of bounds $B$ is a multidimensional box and all procedures presented herein only work with one ReLU at a time. It is possible to devise a procedure for simultaneous tightening of multiple bounds. For example, a procedure could tighten the sum of bounds in the same layer.
\item Strong MILP formulations and a related family of cutting planes for ReLU networks were recently presented in \cite{Anderson2018}. These advancements could extend the applicability of ReLU networks as surrogate models by lowering solution times.
\item Cuts and bound tightening may be used deeper in the MILP solver's B\&B tree. It may be advantageous to employ cheap BT procedures specialized for the structure of ReLU network at selected branches in the branch-and-bound tree.
\item $L_1$ regularization can be utilized during training to encourage sparsity in the parameter matrices, which likely leads to a reduction in MILP solution times \cite{Anderson2018}. Above some level of sparsity, the generalization error will generally start to deteriorate due to the reduction in model capacity. However, as DNNs often are over-parameterized, many connections can typically be dropped without any significant impact on the generalization error.
\item Larger networks may be considered by utilizing specialized solution methods for ReLU networks. In particular, in the literature on robustness and verification of DNNs we find methods such as the modified Simplex algorithm in \cite{Katz2017}, the Lagrangian relaxation-based method in \cite{Krishnamurthy2018}, and methods based on linear and Lipschitz relaxations \cite{Weng2018,Singh2018}. While these works focus on forward bound-propagation in a single network, their ideas may be exploited to efficiently tighten the big-M values in \eqref{eq:nn-milp} or to develop new solution methods for problems with multiple ReLU network surrogate models. 
\end{itemize}
Due to the wide application of ReLU networks for modeling, we believe that it is important to push the practical limitation on network size to allow for more expressive networks as surrogate models in process optimization.


\appendix

%

\section{Forward- vs backward-propagation of bounds}
\label{app:forward-vs-backward-propagation}
Consider $n$ nodes with outputs $x_i$ feeding into a node $y$. The nodes are related as follows
\begin{equation}
y = \sum\limits_{i=1}^{n} w_i x_i + b,
\end{equation}
with parameters $w_i \in \mathbb{R}$ and $b \in \mathbb{R}$. If $x_i$ are outputs of ReLU units we have $l_i \geq 0$ and no sign restrictions on $w_i$. Without loss of generality and to ease the presentation, we here instead require $w_i \geq 0$ and allow negative bounds on $x_i$, which yields an equivalent argument. A forward-propagation of bounds gives $y \in [L_f, U_f]$, with
\begin{equation}
L_f = \sum\limits_{i=1}^{n} w_i l_i + b \quad \text{and } \quad U_f = \sum\limits_{i=1}^{n} w_i u_i + b.
\end{equation}
Thus, it is clear that the bounds on $y$ are linearly dependent on the bounds of the $x_i$-s. A consequence of this is that a tightening of any bound $u_i$ or $l_i$ will have a tightening effect on $U_f$ or $L_f$, respectively.

Now, consider the backward-propagation of bounds on $y$ to node $x_j$. Let $y \leq U_f - \Delta U$, where $\Delta U \geq 0$ is the amount of tightening in the upper bound of $y$. We are interested in how this affects the bound $x_j \leq u_j - \Delta u_j$, where $\Delta u_j \geq 0$ is the change in the upper bound of $u_j$. Using the relationship between the nodes we have that
\begin{equation*}
w_j x_j = y - \sum\limits_{i \neq j} w_i x_i - b
\end{equation*}
The upper bound on $x_j$ is then found to be
\begin{equation*}
\begin{aligned}
w_j (u_j - \Delta u_j) &= U_f - \Delta U - \sum\limits_{i \neq j} w_i l_i - b \\
&= \sum\limits_{i=1}^{n} w_i u_i + b - \Delta U - \sum\limits_{i \neq j} w_i l_i - b \\
&= \sum\limits_{i=1}^{n} w_i u_i - \Delta U - \sum\limits_{i \neq j} w_i l_i 
\end{aligned}
\end{equation*}
Rearranging the last equation yields
\begin{equation*}
\Delta u_j = \frac{\Delta U}{w_j} - \sum\limits_{i \neq j} \frac{w_i}{w_j} (u_i - l_i)
\end{equation*}
Tightening of the upper bound on $x_j$ occurs only when $\Delta u_j \geq 0$, and we obtain the requirement
\begin{align}
& \frac{\Delta U}{w_j} - \sum\limits_{i \neq j} \frac{w_i}{w_j} (u_i - l_i) \geq 0 \notag \\
\implies & \Delta U \geq \sum\limits_{i \neq j} w_i (u_i - l_i) \geq 0
\label{eq:backward-prop-requirement-upper-bound}
\end{align}

An equivalent argument can be made for the lower bound $x_j \geq l_j + \Delta l_j$, given a tightening $y \geq L_f + \Delta L$. Requiring $\Delta l_j \geq 0$ and $\Delta L \geq 0$ gives 
\begin{equation}
\Delta L \geq \sum\limits_{i \neq j} w_i (u_i - l_i) \geq 0.
\label{eq:backward-prop-requirement-lower-bound}
\end{equation}

Inequalities \eqref{eq:backward-prop-requirement-lower-bound} and \eqref{eq:backward-prop-requirement-upper-bound} tell us that a considerable change in the bound of $y$ is required for it to have any tightening effect on the bounds of preceding nodes. To better understand these requirements we may cast them as a relative change in the bounds by considering the quantity $\Delta / (U_f - L_f)$, where $\Delta$ is the change in the bound of $y$ (either $\Delta U$, $\Delta L$, or a combination of these). Since $U_f - L_f \geq 0$, we may write
\begin{align*}
\frac{\Delta}{U_f - L_f} &\geq \frac{1}{U_f - L_f} \sum\limits_{i \neq j} w_i (u_i - l_i) \\
&= \frac{\sum\limits_{i \neq j} w_i (u_i - l_i)}{\sum\limits_{i =1}^{n} w_i (u_i - l_i)} \\
&= 1 - \frac{w_j (u_j - l_j)}{\sum\limits_{i =1}^{n} w_i (u_i - l_i)} \\
&= 1 - \delta,
\end{align*}
where we introduced $\delta \in [0, 1]$ in the last line. It now becomes evident that a relatively large change $\Delta$ is required for backward-propagation of bounds to have any effect. We see that $\delta \approx 1$ only if $w_j \gg w_i$ or $u_j - l_j \gg u_i - l_i$ for all $i \neq j$, meaning that the relationship between $x_j$ and $y$ dominates the relationship between $y$ and the other nodes.

We end this discussion by considering the special case where $w_i = w_j$, $u_i = u_j$, and $l_i = l_j$  for all $i, j$. In this case $\delta$ simplifies to $\delta = 1/n$ and we get
\begin{equation*}
\frac{\Delta}{U_f - L_f} \geq 1 - 1/n.
\end{equation*}
From this expression we see that for wide layers (large $n$) with ``balanced'' nodes, the tightening $\Delta$ must be very large for it to have any effect on preceding nodes. E.g., for $n=100$, a reduction of $99\%$ is required in the bounds on $y$.

\section{Statistics for measuring bound tightness}
\label{app:bound-statistics}
We introduce two statistics to measure the quality of the bounds produced by a bound tightening procedure. We measure the mean absolute distance (MAD) of a set of bounds $B$ as
\begin{equation}
MAD(B) := \sum\limits_{k=0}^{K} n_k^{-1} \sum\limits_{j=1}^{n_k} |u_{j}^{k} - l_{j}^{k}|,
\end{equation}
where $u_{j}^{k}$ and $l_{j}^{k}$ are the upper and lower bound on $x_j^k - s_j^k$ for node $(j,k)$.

To compare bounds, we use the mean relative distance (MRD), which measures the relative tightness of bounds $B$ compared to the best bounds $B^\star$ and the suboptimal LRR bounds denoted $B^-$. The MRD is defined as follows:
\begin{equation}
MRD(B, B^{\star}, B^{-}) := 100 \frac{|\text{MAD}(B) - \text{MAD}(B^\star)|}{|\text{MAD}(B^-) - \text{MAD}(B^\star)|}.
\end{equation}
For a set of optimal bounds $B = B^*$, $MRD(B, B^{\star}, B^{-}) = 0$. While for bounds $B = B^-$, the MRD is 100. The bound tightening procedures used in this study are initialized with $B^0$, and thus always obtain an MRD in the range $[0, 100]$. Note that we define the MRD to be zero whenever the denominator is zero.

\section{Bound tightening results}
\label{app:bound-tightening-results}


\begin{table}[ht!]
\footnotesize
\caption{Bound tightening results for \ref{eq:quad-problem}. The reported values are the average results over ten models, for $n=1,\ldots,6$.}
\begin{tabularx}{\textwidth}{ll*4{>{\raggedleft\arraybackslash}X}}
\toprule
n & BT procedure & Time (s) & Dead neurons (\%) & MAD & MRD (\%) \\
\midrule 
1 & LRR & 0.013 & 29.3 & 2.075 & 100.0 \\  
  & RR & 0.026 & 41.3 & 1.254 & 7.4 \\  
  & LR(1) & 0.054 & 32.7 & 1.695 & 39.2 \\  
  & SEMI-LR(1) & 0.138 & 42.0 & 1.238 & 4.5 \\  
  & NO-R(1) & 0.070 & 43.3 & 1.179 & 0.0 \\  
  \noalign{\smallskip\smallskip}
2 & LRR & 0.024 & 9.3 & 5.832 & 100.0 \\  
  & RR & 0.078 & 12.0 & 2.916 & 21.3 \\  
  & LR(1) & 0.374 & 12.3 & 2.749 & 15.2 \\  
  & SEMI-LR(1) & 1.575 & 12.7 & 2.598 & 10.7 \\  
  & NO-R(1) & 1.959 & 16.7 & 2.228 & 0.0 \\ 
  \noalign{\smallskip\smallskip}
3 & LRR & 0.065 & 9.0 & 9.879 & 100.0 \\  
  & RR & 0.293 & 9.2 & 4.865 & 27.6 \\  
  & LR(1) & 3.772 & 10.2 & 3.866 & 12.4 \\  
  & SEMI-LR(1) & 17.42 & 10.5 & 3.747 & 9.6 \\  
  & NO-R(1) & 36.65 & 13.8 & 3.245 & 0.0 \\  
  \noalign{\smallskip\smallskip}  
4 & LRR & 0.082 & 4.9 & 14.76 & 100.0 \\  
  & RR & 0.403 & 4.9 & 6.53 & 25.1 \\  
  & LR(1) & 7.717 & 6.1 & 4.75 & 7.6 \\  
  & SEMI-LR(1) & 32.02 & 6.1 & 4.68 & 6.8 \\  
  & NO-R(1) & 84.29 & 8.9 & 4.13 & 0.0 \\  
  \noalign{\smallskip\smallskip}  
5 & LRR & 0.164 & 4.0 & 67.47 & 100.0 \\  
  & RR & 2.252 & 4.1 & 18.15 & 17.4 \\  
  & LR(1) & 95.59 & 7.5 & 8.79 & 1.6 \\  
  & SEMI-LR(1) & 181.2 & 6.9 & 9.18 & 2.2 \\  
  & NO-R(1) & 303.9 & 7.0 & 8.55 & 1.1 \\  
  \noalign{\smallskip\smallskip}
6 & LRR & 0.286 & 3.4 & 48.18 & 100.0 \\  
  & RR & 5.970 & 3.4 & 13.41 & 15.4 \\  
  & LR(1) & 195.8 & 4.2 & 7.14 & 0.0 \\  
  & SEMI-LR(1) & 335.0 & 3.4 & 9.59 & 6.1 \\  
  & NO-R(1) & 472.9 & 3.5 & 8.49 & 3.4 \\  
\bottomrule \noalign{\smallskip}
\end{tabularx}
\label{tab:bt-results-quadratic}
\end{table}


\begin{table}[ht!]
\footnotesize
\caption{Bound tightening results for shallow well networks with layers (1, 20, 20, 1). The reported values are the average results of the eight models.}
\begin{tabularx}{\textwidth}{l*5{>{\raggedleft\arraybackslash}X}}
\toprule
BT method & Time (s) & Dead neurons (\%) & MAD & MRD (\%) \\
\midrule 
LRR & 0.020 & 0.0 & 0.2192 & 100.0 \\  
RR & 0.062 & 0.7 & 0.2027 & 12.0 \\  
LR(60) & 0.097 & 1.8 & 0.2031 & 14.6 \\  
SEMI-RR(60) & 0.548 & 1.8 & 0.2027 & 12.0 \\  
NO-R(60) & 0.104 & 1.8 & 0.2014 & 0.0 \\  
NO-R & 0.104 & 1.8 & 0.2014 & 0.0 \\  
\bottomrule \noalign{\smallskip}
\end{tabularx}
\label{tab:bt-results-shallow-well-nets}
\end{table}

\begin{table}[ht]
\footnotesize
\caption{Bound tightening results for deep well networks with layers (1, 10, 10, 10, 10, 1). The reported values are the average results of the eight models.}
\begin{tabularx}{\textwidth}{l*5{>{\raggedleft\arraybackslash}X}}
\toprule
BT method & Time (s) & Dead neurons (\%) & MAD & MRD (\%) \\
\midrule 
LRR & 0.022 & 0 & 0.5723 & 100.0 \\
RR & 0.065 & 0 & 0.5457 & 11.5 \\
LR(60) & 0.183 & 0 & 0.5461 & 12.5 \\
SEMI-RR(60) & 0.503 & 0 & 0.5457 & 11.5 \\
NO-R(60) & 0.156 & 0 & 0.5421 & 0.0 \\
NO-R & 0.156 & 0 & 0.5421 & 0.0 \\
\bottomrule \noalign{\smallskip}
\end{tabularx}
\label{tab:bt-results-deep-well-nets}
\end{table}

\begin{table}[ht]
\footnotesize
\caption{Bound tightening results for shallow riser networks with layers (4, 50, 50, 1). The reported values are the average results of the two models.}
\begin{tabularx}{\textwidth}{l*5{>{\raggedleft\arraybackslash}X}}
\toprule
BT method & Time (s) & Dead neurons (\%) & MAD & MRD (\%) \\
\midrule 
LRR & 0.093 & 4.0 & 2.915 & 100.0 \\
RR & 0.554 & 7.0 & 1.998 & 46.0 \\
LR(60) & 8.889 & 16.0 & 1.571 & 27.2 \\
SEMI-RR(60) & 55.29 & 16.0 & 1.521 & 23.6 \\
NO-R(60) & 306.0 & 19.5 & 1.122 & 0.0 \\
NO-R & 306.0 & 19.5 & 1.122 & 0.0 \\
\bottomrule \noalign{\smallskip}
\end{tabularx}
\label{tab:bt-results-shallow-riser-nets}
\end{table}

\begin{table}[ht]
\footnotesize
\caption{Bound tightening results for deep riser networks with layers (4, 20, 20, 20, 20, 20, 1). The reported values are the average results of the two models.}
\begin{tabularx}{\textwidth}{l*5{>{\raggedleft\arraybackslash}X}}
\toprule
BT method & Time (s) & Dead neurons (\%) & MAD & MRD (\%) \\
\midrule 
LRR & 0.061 & 0 & 13.418 & 100.0 \\
RR & 0.991 & 1 & 4.120 & 27.4 \\
LR(60) & 109.3 & 4 & 2.449 & 14.7 \\
SEMI-RR(60) & 134.9 & 4 & 2.377 & 14.2 \\
NO-R(60) & 1771 & 11 & 0.678 & 0.1 \\
NO-R & 5080 & 11 & 0.662 & 0.0 \\
\bottomrule \noalign{\smallskip}
\end{tabularx}
\label{tab:bt-results-deep-riser-nets}
\end{table}

\clearpage
\section*{References}
\bibliography{references}

\end{document}